\documentclass[12 pt]{article}
\usepackage{amssymb,amsmath,amsfonts,amsthm}
\usepackage{xcolor}
\newcommand\BBP{{\mathbb {P}}}

\newcommand\E{{\mathbb {E}}}

\newtheorem {Lemma}{Lemma}[section]
\newtheorem {Theorem}{Theorem}[section]

\theoremstyle{definition}

\newtheorem{Remark}{Remark}[section]

\newtheorem{Comment}{Comment}[section]

\newcommand\beq{\begin{equation}}
\newcommand\eeq{\end{equation}}

\begin{document}

\title{Berry-Esseen type bounds for the Left Random Walk on $GL_d({\mathbb R})$ under polynomial moment conditions}

\author{C. Cuny\footnote{Christophe Cuny, Univ Brest,  UMR CNRS 6205, LMBA, 6 avenue Victor Le Gorgeu, 29238 Brest}, J. Dedecker\footnote{J\'er\^ome Dedecker, Universit\'e  Paris Cit\'e, CNRS, MAP5, UMR 8145,
45 rue des  Saints-P\`eres,
F-75006 Paris, France.}, F. Merlev\`ede \footnote{Florence Merlev\`ede, LAMA,  Univ Gustave Eiffel, Univ Paris Est Cr\'eteil, UMR 8050 CNRS,  \  F-77454 Marne-La-Vall\'ee, France.}
and 
M. Peligrad \footnote{Magda Peligrad, Department of Mathematical Sciences, University of Cincinnati, PO Box 210025, Cincinnati, Oh 45221-0025, USA.}}

\maketitle

\begin{abstract} Let $A_n= \varepsilon_n \cdots \varepsilon_1$,  where $(\varepsilon_n)_{n \geq 1}$ is a sequence of  independent random matrices taking values in $ GL_d(\mathbb R)$, $d \geq 2$, with  common distribution $\mu$. In this paper, under standard assumptions on $\mu$ (strong irreducibility and proximality), we prove Berry-Esseen type theorems for $\log ( \Vert A_n \Vert)$ when $\mu$ has a polynomial moment. More precisely, we get the rate  $ ( (\log n ) / n)^{q/2-1}$ when  $\mu$ has a moment of order $q \in ]2,3]$ and the  rate $1/ \sqrt{n} $ when  $\mu$ has a moment of order $4$, which significantly improves  earlier results in this setting. 
\end{abstract}



{

\section{Introduction}

\setcounter{equation}{0}

Let $(\varepsilon_n)_{n \geq 1}$ be independent random matrices taking values in $G= GL_d(\mathbb R)$, $d \geq 2$ (the group of invertible $d$-dimensional real matrices) with common distribution $\mu$. Let $\Vert \cdot \Vert$ be the euclidean norm on ${\mathbb R}^d$, and for every $A \in GL_d(\mathbb R)$, let $\|A\|=\sup_{x, \|x\|=1} \|A x \|$. We shall say that $\mu $ has a moment of order $p \geq 1$ if
\[
\int_G (\log N(g) )^p d \mu(g) < \infty \, , 
\]
where $N(g) := \max ( \Vert g \Vert , \Vert g^{-1} \Vert)$.  

Let $A_n= \varepsilon_n \cdots \varepsilon_1$. It follows from Furstenberg and Kesten \cite{FK} that, if $\mu$ admits a moment of order $1$ then 
\beq \label{SL1}
\lim_{n \rightarrow \infty} \frac{1}{n}  \log \Vert A_n \Vert = \lambda_{\mu} \, \text{ ${\mathbb P}$-a.s.},
\eeq
where $ \lambda_{\mu} := \lim_{n \rightarrow \infty} n^{-1}  \E \log \Vert A_n  \Vert $ is the so-called first Lyapunov exponent.

Let now $X:= P({\mathbb R}^d)$  be the projective space of ${\mathbb R}^d $ and write ${\bar x}$ as the projection of $x \in {\mathbb R}^d -\{0\}$ to $X$.  An element $A$ of $G= GL_d(\mathbb R)$ acts on  the projective space $X$ as follows: $A \bar x = \overline{Ax}$.  Let $\Gamma_\mu$ be the closed semi-group generated by the support of $\mu$. We say that $\mu$ is proximal if $\Gamma_\mu$ contains a matrix  that admits  a unique (with multiplicity $1$) eigenvalue of maximal modulus. We say that $\mu$ is strongly irreducible if no proper union of subspaces of ${\mathbb R}^d$ is invariant by $\Gamma_\mu$. Throughout the paper, we  assume that $\mu$ is strongly irreducible and proximal.
% (see for instance \cite{BL} or \cite{BQ} for a clear exposition of these notions), 
In particular,  there exists a unique invariant measure $\nu$ on ${\mathcal B} (X)$ with respect to $ \mu$, meaning that for any continuous and bounded function $h$ from $X$ to $\mathbb R$,
\beq \label{defnu}
\int_X h(x) d \nu(x) = \int_G \int_X h( g \cdot x ) d \mu(g) d \nu(x) \, .
\eeq
Note that, since $\mu$ is assumed to be strongly irreducible, the following strong law holds (see for instance \cite{BL}, Proposition 7.2 page 72): for any $x \in {\mathbb R}^d -\{0\}$, 
\beq \label{SL2}
\lim_{n \rightarrow \infty} \frac{1}{n}  \log \Vert  A_n  x \Vert = \lambda_{\mu} \, \text{ ${\mathbb P}$-a.s.}
\eeq
%The left random walk of law $\mu$ is the process defined by $W_0:=\varepsilon_0$ and $W_n %=\varepsilon_n W_{n-1}$ for $n \geq 1$ where  $\varepsilon_0$ is a random variable with values in 
%$X$ and independent of the sequence  $(\varepsilon_n)_{n \geq 1}$. 
To specify the rate of convergence in the laws of large numbers \eqref{SL1} and \eqref{SL2}, it is then natural to address the question of the Central Limit Theorem for the two sequences 
$\log \| A_n \| - n\lambda_\mu$ and $\log \| A_n  x \| - n\lambda_\mu$.  To specify the limiting variance in these central limit theorems, let us introduce some notations: $W_0$ will denote a random variable with values in the projective space $X$, independent of $(\varepsilon_n)_{n \geq 1}$  and with distribution $\nu$. By the invariance of $\nu$, we see that the process $(A_n W_0)_{n \geq 1}$ is a strictly stationary process. Denote also by 
$V_0$ a random variable such that $\Vert V_0 \Vert =1 $ and ${\bar V}_0 = W_0$. Setting, $S_n = \log \| A_nV_0 \| - n\lambda_\mu$,  Benoist and Quint \cite{BQ}  proved that if $\mu$ has a moment of order $2$, then  
\beq \label{definitions2}
\lim_{n \rightarrow \infty} \frac{1}{n} \E ( S_n^2) = s^2 >0 \, , 
\eeq
\beq \label{BQCLT}
 \lim_{n \rightarrow \infty} \sup_{t \in {\mathbb R}} \sup_{x, \|x\|=1}\left | {\mathbb P}  \left ( 
 \log \| A_n   x \| -  n \lambda_\mu \leq t  \sqrt n  \right ) - \Phi (t/ s)  \right | =0  \, ,
\eeq
and
\beq \label{BQCLT2}
 \lim_{n \rightarrow \infty} \sup_{t \in {\mathbb R}} \left | {\mathbb P}  \left ( 
 \log \| A_n    \| -  n \lambda_\mu \leq t  \sqrt n  \right ) - \Phi (t/ s)  \right | =0  \, ,
\eeq
where $\Phi$ is the cumulative distribution function of a standard normal distribution.  Let us mention that   \eqref{BQCLT}  has been firstly established by Le Page \cite{LP} under an exponential moment for $\mu$ (meaning that  $\int_G (N(g))^\alpha  d \mu(g) < \infty$ for some $\alpha>0$, see also  \cite{GR}) and   then  by Jan \cite{Jan}  under the condition that  $\mu$ has a moment of order $p>2$.

In the present paper, we are interested in Berry-Esseen type bounds in these central limit theorems, under polynomial moments for $\mu$ (more precisely we shall focus on the case of moments of order $q \in ]2,3]$ or $q=4$). Before giving our main results, let us briefly describe the previous works on this subject. 

When $\mu$ has an exponential moment, Le Page \cite{LP} proved the following inequality: there exists a positive constant $C$ such that 
\beq \label{Lepage}
\sup_{t \in {\mathbb R}} \sup_{x, \|x\|=1}\left | {\mathbb P}  \left ( 
 \log \| A_n   x \| -  n \lambda_\mu \leq t  \sqrt n  \right ) - \Phi (t/ s)  \right | \leq  C v_n \,  \text{ with } \, v_n= \frac{1}{\sqrt{n}} \, .
\eeq
Still in the case of exponential moments, Edgeworth expansions  (a strengthening of the Berry-Esseen theorem) have been recently obtained by Fernando and P\`ene \cite{FP20} and Xiao et al. \cite{XGLeuropean}.  In these three last   papers, the assumption that $\mu$ has an exponential moment is crucial since it  allows to use the strength of the so-called  Nagaev-Guivarc'h perturbation method. Indeed, in case of exponential moments, the associated complex perturbed transfer operator has  spectral gap properties.

Now,  under the assumption that all the moments of order $p$ of $\mu$ are finite,  Jan \cite{Jan} obtained the rate  $v_n= n^{-1/2+ \varepsilon}$   for any $\varepsilon >0$ in \eqref{Lepage}. Next, 
Cuny et al. \cite{CDJ} gave an upper bound of order $v_n =  n^{-1/4} \sqrt{\log n}$ in \eqref{Lepage} provided $\mu$ has a moment of order $3$ (as a consequence of an upper bound  of order 
$ n^{-1/2} \log n$ for the Kantorovich metric). More recently, Jirak \cite{Ji20} proved that, if $\mu$ has a moment of order $p>8$, then
there exists a positive constant $C$ such that 
\beq \label{Jirak}
\sup_{t \in {\mathbb R}} \left | {\mathbb P}  \left ( 
\log \| A_n  V_0 \| - n \lambda_\mu \leq t  \sqrt n \right ) - \Phi (t/ s)  \right | \leq  C v_n \,  \text{ with } \, v_n= \frac{1}{\sqrt{n}}  \, .
\eeq
This result is based on some refinements of the arguments developed in a previous paper of the same author (see \cite{Ji}), and then on a completely different method than   the   perturbation method for the transfer operator. 
Since our proofs will use a similar scheme let us briefly explain it. First, due to the cocycle property (see the beginning of Section \ref{BEbounds}), $\log \| A_n  V_0 \| - n \lambda$ is written as a partial sum associated with functions of a stationary Markov chain, which can be viewed also as   a function of iid random elements (see also \cite{CDM}).  Using the conditional expectation, the underlying random variables are then approximated  by $m$-dependent variables, say $X_{k,m}$. Next, to break the dependence, a blocking procedure is used and  the partial sum $\sum_{k=1}^n X_{k,m}$ is decomposed into two terms. The first one can be rewritten as the  sum of random variables  which are defined as  blocks, say $Y_j^{(1)}$,  of size $2m$ of the $X_{k,m}$'s. These  random blocks have the following property: conditionally to ${\mathbb F}_m $ (a particular $\sigma$-algebra generated by a part of the $\varepsilon_k$'s), they are independent. In addition, for any bounded measurable function $h$, the random variables $Z_j= \E (h ( Y_j^{(1)})  | {\mathbb F}_m) $  are one-dependent. On another hand,  the second term in the decomposition of $\sum_{k=1}^n X_{k,m}$ is ${\mathbb F}_m $-measurable and can be written as a sum of independent blocks of the initial random variables.  For both terms in the decomposition, the conditional  independence  of the blocks comes from the independence of the $\varepsilon_k$'s. The next steps of the proof consist first of all in working conditionally   to ${\mathbb F}_m $  and  then in giving  suitable upper bounds for the conditional  characteristic function of the blocks $Y_j^{(1)}$.

%where $V_0$ is independent of $(\varepsilon_n)_{n \geq 1}$ and   such that $\Vert V_0 \Vert =1$ and ${\overline {V_0}}$ is  distributed according to the invariant distribution $\nu$ on $X$. 

Concerning matrix norms, we first note that the Berry-Esseen bound of order $n^{-1/4} \sqrt{\log n}$ under a moment of order $3$ is still valid for $\log \| A_n \| -n  \lambda_\mu$ instead of $\log \| A_n  x \| - n \lambda_\mu$ (see the discussion in Section 8 of \cite{CDJ}). Moreover, if $\mu$ has an exponential moment, Xiao et al. \cite{XGL} proved that there exists a positive constant $C$ such that	
\beq \label{Liu}
\sup_{t \in {\mathbb R}} \left | {\mathbb P}  \left ( 
\log \| A_n \| - n \lambda_\mu \leq t  \sqrt n \right ) - \Phi (t/ s)  \right | \leq  C w_n \,  \text{ with } \,w_n =   \frac{ \log n}{\sqrt{n}} \, .
\eeq
Note that in \cite{XGL}, the authors also  proved a similar upper bound for $\log (\rho (A_n) )$ where  $\rho (A_n)$ is the spectral radius of $A_n$. 

\smallskip

In the present paper, we prove that: 
\begin{itemize}
\item If $\mu$ has a moment of order $q \in ]2,3]$, then the rate in \eqref{Lepage}  (and then in \eqref{Jirak}) is $ v_n = ( \log n/n)^{q/2-1}$ and the rate in 
\eqref{Liu} is $w_n =( \log n/n)^{q/2-1}$.
\item If $\mu$ has a moment of order 4, then the rate in \eqref{Lepage} (and then in \eqref{Jirak})  is $ v_n=  n^{-1/2} $ and the rate in 
\eqref{Liu} is $w_n = n^{-1/2} $.
\end{itemize}

To prove these results, we follow the blocking approach used in Jirak \cite{Ji,Ji20} (and described above), but with substantial changes. We refer to Comment \ref{comment31}  to have a flavor  of  them. One of the main changes is the use of the dependency coefficients defined  in 
\cite{CDJ} (see also  \eqref{delta-def} in Section \ref{sectionproofs}) which are well adapted to the study of the process $( \log \| A_n   x \| -  n \lambda_\mu)_{n \geq 1}$, instead of the coupling coefficients used in \cite{Ji20}.

\smallskip

The paper is organized as follows. In Section \ref{BEbounds}, we state our main results about  Berry-Esseen type bounds in the context of left random walks when $\mu$ has  either a moment of order $q \in ]2,3]$ or a moment of order $4$. All the proofs are  postponed to Section \ref{sectionproofs}. Some technical lemmas used in the proofs are stated and proved in Section \ref{TL}. 

\smallskip

In the rest of the paper, we shall use the following notations:  for two sequences $(a_n)_{n \geq 1}$ and $(b_n)_{n \geq 1}$ of positive reals, $a_n \ll b_n$ means that there exists a positive constant $C$ not depending on $n$ such that $a_n \leq C b_n$ for any $n\geq 1$. Moreover, given a $\sigma$-algebra ${\mathcal F}$, we shall often use the notation $\E_{{\mathcal F}} ( \cdot) = \E ( \cdot |{\mathcal F}) $. 

\begin{Remark}
After this article was submitted, we became aware of the paper
by Dinh, Kaufmann and Wu \cite{DKW}, in which the authors obtain the bound \eqref{Lepage} with $v_n = n^{-1/2}$ 
when $\mu$ has a moment of order $3$, but only in the case $d = 2$. Note that, in the same paper and still
in the case $d = 2$, a Local Limit Theorem is also established for $\log  \Vert A_n x \Vert$. 
\end{Remark}

\section{Berry-Esseen bounds} \label{BEbounds}

\setcounter{equation}{0}

Recall the notations in the introduction: let $(\varepsilon_n)_{n \geq 1}$ be independent random matrices taking values in $G= GL_d(\mathbb R)$, $d \geq 2$, with common distribution $\mu$. Let $A_n= \varepsilon_n \cdots \varepsilon_1$  for $n\geq 1$, and $A_0=$Id.  We assume that $\mu$ is strongly irreducible and proximal, and we denote by $\nu$ the unique  distribution on $X=P({\mathbb R}^d)$ satisfying \eqref{defnu}. 

Let now $V_0$ be a random variable independent of $(\varepsilon_n)_{n \geq 1}$, taking values in ${\mathbb R}^d$, such that $\Vert V_0 \Vert =1$ and ${\overline {V_0}}$ is distributed according to  $\nu$.

The behavior of $\log \| A_n  V_0 \| - n \lambda_\mu $ (where $\lambda_\mu$ is the first Lyapunov exponent defined right after \eqref{SL1})  can be handled with the help of an additive cocycle,  which can also be viewed as a function of a stationary Markov chain. More precisely, let $W_0=\overline {V_0}$ (so that $W_0$ is distributed according to $\nu$), and let $W_n = \varepsilon_n W_{n-1}=A_n W_0$ for  any integer $n \geq 1$. By definition of $\nu$, the sequence $(W_n)_{n \geq 0}$ is a strictly stationary Markov chain with values in $X$.  Let now, for any integer $k \geq 1$, 
\beq  \label{defMCXk}
X_k := \sigma (\varepsilon_k, W_{k-1} ) - \lambda_{\mu} = \sigma (\varepsilon_k, A_{k-1} W_0 ) - \lambda_{\mu} \, , 
\eeq 
where,   for any $g \in G$ and any ${\bar x} \in X$,
\[
\sigma( g , {\bar x} ) = \log \Big ( \frac{\Vert g \cdot x \Vert }{ \Vert x \Vert }\Big ) \, .
\]
Note that $\sigma$ is an additive cocycle in the sense that $\sigma ( g_1 g_2,   {\bar x})  = \sigma ( g_1,    g_2 {\bar x})  +  \sigma ( g_2,  {\bar x}) $.  Consequently
% Hence, if we set $A_0={\rm Id}$ and for every $n \geq 1$, $A_n = \varepsilon_n \cdots 
%\varepsilon_1$, we get $X_k = \sigma (\varepsilon_k, A_{k-1}W_0) - \lambda_{\mu}$. So that 
\beq  \label{defofSn}
S_n = \sum_{k=1}^n X_k = \log \Vert  A_n V_0 \Vert  - n \lambda_{\mu}\, .
\eeq

With the above notations, the following Berry-Esseen bounds hold. 
\begin{Theorem} \label{thmq=3}
Let $\mu$ be a proximal and strongly irreducible probability measure on ${\mathcal B} (G)$. Assume that $\mu$ has a finite moment of order $q \in ]2,3]$. Then $n^{-1} \E (S_n^2) \rightarrow s^2>0$ as $n \rightarrow \infty$ and, setting $ \displaystyle v_n =  \Big (  \frac{\log n}{n }  \Big )^{q/2-1}$, we have 
\beq \label{ineBE1}
\sup_{y \in {\mathbb R}} \Big | {\mathbb P} \big ( S_n  \leq y \sqrt{n} \big ) - \Phi (y/ s)  \Big | \ll v_n \, , \eeq
\beq \label{ineBE1bis}
  \sup_{y \in {\mathbb R}} \Big | {\mathbb P} \big ( \log (\Vert A_n \Vert ) -  n \lambda_\mu  \leq y \sqrt{n} \big ) - \Phi (y/ s)  \Big | \ll v_n   \, , 
\eeq
and
\beq \label{ineBE1ter}
 \sup_{x , \Vert x \Vert =1}\sup_{y \in {\mathbb R}} \Big | {\mathbb P} \big ( \log \Vert  A_n x \Vert  - n \lambda_{\mu}  \leq y \sqrt{n} \big ) - \Phi (y/ s)  \Big | \ll v_n \, . \eeq
\end{Theorem}

\begin{Remark} \label{remonvar}
As mentioned in the introduction, the fact that  
$n^{-1} \E (S_n^2) \rightarrow s^2 >0$ has been proved by Benoist and Quint \cite{BQ} (see Item $(c)$ of  their Theorem 4.11). Let us mention that we also have 
$s^2 =\E (X_1^2) + 2 \sum_{k \geq 2} \E  (X_1X_k) $, which follows for instance from the proof of item $(ii)$ of Theorem 1 in \cite{CDJ}. 
\end{Remark}

\begin{Remark} \label{CRAS}
The results of Theorem \ref{thmq=3} are used in the article \cite{CDMP} to obtain Berry-Esseen type bounds for the matrix coefficients and for the spectral radius, that is for the quantities 
\begin{align*}
  \sup_{ \Vert x \Vert=  \Vert y \Vert=1}\sup_{t \in {\mathbb R}} \Big | {\mathbb P} \big ( \log |\langle A_n x, y\rangle  |-  n \lambda_\mu  \leq t \sqrt{n} \big ) - \Phi (t/ s)  \Big |  \, ,\\
 \text{and} \quad  \sup_{t \in {\mathbb R}} \Big | {\mathbb P} \big ( \log \lambda_1(A_n)  -  n \lambda_\mu  \leq t \sqrt{n} \big ) - \Phi (t/ s)  \Big | \, ,
\end{align*} 
where $\lambda_1(A_n)$ is the greatest modulus of the eigenvalues of the matrix $A_n$. In \cite{CDMP}, only the case of polynomial moments of order $q\geq 3$ is considered, but it is actually possible to obtain bounds for moments $q>2$ thanks to Theorem \ref{thmq=3}.  More precisely, for both quantities, the rates are 
\begin{itemize}
\item $v_n=(\log n /n)^{q/2-1}$ if $\mu$ has a finite moment of order $q \in ]2, (3+\sqrt 5)/2]$;
\item $v_n=1/n^{(q-1)/2q}$ if $\mu$ has a finite moment of order $q> (3+\sqrt 5)/2$.
\end{itemize}
\end{Remark}

\medskip

Now if  $\mu$ has a finite moment of order $4$ then the following result holds: 

\begin{Theorem}  \label{thmq=4}
Let $\mu$ be a proximal and strongly irreducible probability measure on ${\mathcal B} (G)$. Assume that $\mu$ has a finite moment of order $4$. Then $n^{-1} \E (S_n^2) \rightarrow s^2>0$ as $n \rightarrow \infty$ and \eqref{ineBE1}, 
\eqref{ineBE1bis} and \eqref{ineBE1ter} hold with $v_n= 1/\sqrt{n}$. 
\end{Theorem}

Recall that the classical Berry-Esseen theorem for independent random variables, which corresponds to the case $d=1$ in our setting, provides the rate 
 $1/{\sqrt n}$ under a finite moment of order $3$. For $q=3$, Thorem \ref{thmq=3} provides the rate $\sqrt{(\log n) /n}$, so  one may wonder whether the conclusion of Theorem \ref{thmq=4} holds when $\mu$ has a moment of order $3$ only. 
 
 Note also that we have chosen to focus on the cases where $\mu$ has a finite moment of order $q \in ]2,3]$  (since it corresponds to the usual moment assumptions for the Berry-Esseen theorem in the iid case)  or a finite moment of order $4$ (since in this case we reach the rate $1/\sqrt{n}$), but we infer from the proofs that if $\mu$ has a finite moment of order $q \in ]3,4[$ then the above results hold with $v_n=  ( \log n  )^{(4-q)/2} /{\sqrt{n}} $. 

\section{Proofs} \label{sectionproofs}

\setcounter{equation}{0}

\subsection{Proof of Theorem \ref{thmq=3}} 

As usual, we shall denote by $X_{k,{\bar x}}$ the random variable $X_k$  defined by \eqref{defMCXk} when the Markov chain $(W_n)_{n \geq 0}$ starts from $\bar x\in X$. We then define $S_{n,{\bar x}}:= \log \big(\|A_n x\|/\|x\|\big) -n \lambda_\mu= \sum_{k=1}^n X_{k,{\bar x}}$.  We shall first prove the upper bound \eqref{ineBE1} in Section \ref{subsection1} and then the upper bounds \eqref{ineBE1bis} and \eqref{ineBE1ter} in Sections 
\ref{subsection2} and \ref{subsection3} respectively.

\subsubsection{Proof of the upper bound  \eqref{ineBE1}} \label{subsection1}

 $ $
 
 \smallskip
 
As usual, the proof is based on  the so-called Berry-Esseen smoothing inequality (see e.g. \cite[Ineq. (3.13)  p. 538]{Fe71}) stating that,  
there exists $C>0$ such that for any positive $T$ and any integer $n\ge 1$, 
\begin{equation} \label{IneBE}
\sup_{x \in {\mathbb R}} \Big | {\mathbb P}  \big ( S_{n}  \leq x \sqrt{n} \big ) - \Phi (x/ s)  \Big |  \le C\Big(  \int_{-T}^T \frac{ \big |  \E \big (  {\rm e}^{{\rm i} \xi  S_n /{\sqrt n}  }\big ) - {\rm e}^{-\xi^2s^2/2}  \big |}{|\xi|} {\rm d}\xi  + T^{-1}\, \Big)\, , 
\end{equation}
 where we recall that $S_n$ has been defined in  \eqref{defofSn}.  
 
 To take care of the characteristic function of $S_n /{\sqrt n} $ we shall take advantage of the fact  that $X_k $ is a function of a stationary Markov chain generated by the iid random elements $(\varepsilon_i)_{i \geq 1}$. As in \cite{Ji}, the first steps of the proof consist in approximating  the $X_k$'s by $m$-dependent random variables $X_{k,m}$, and then in suitably decomposing the partial sum associated with the $X_{k,m}$'s. This is the subject of the following paragraph.  

\smallskip

\noindent {\it Step 0. Notations and Preliminaries.}   We shall   adopt most of the time  the same notations as in Jirak \cite{Ji}. Let ${\mathcal E}_{i}^j = \sigma (\varepsilon_i, \ldots, \varepsilon_j)$ for $i \leq j$, and $m$ be a positive integer that will be specified later. For any $k \geq m$, let 
\beq  \label{defXkm}
X_{k,m} = \E ( X_k| {\mathcal E}_{k-m+1}^k ) := f_m ( \varepsilon_{k-m+1}, \ldots, \varepsilon_k) \, , 
\eeq
where $f_m$ is  a measurable function.  More precisely, we have
\[
X_{k,m} =  \int_X \sigma ( \varepsilon_k, A_{k-1}^{k-m+1} {\bar x} ) d \nu ({\bar x})   - \lambda_{\mu} \, , 
\]
where we used the notation $A_{j}^i = \varepsilon_j \cdots \varepsilon_i$ for $i \leq j$.  Note that $\E(X_{k,m}) =0$.

Next, let $N$ be the positive integer such that  $n = 2 N m + m' $ with $0 \leq m' \leq 2m-1$. 
The integers $N$ and $m$  are such that $N \sim \kappa_1 \log n$ (where $\kappa_1$ is a positive constant specified later) and   $m \sim ( 2\kappa_1)^{-1} n (\log  n)^{-1}$ (see \eqref{selectionofN} for the selection of $\kappa_1$). 
Define now the following $\sigma$-algebra
\beq  \label{defbbFm}
{\mathbb F}_m = \sigma (  (  \varepsilon_{(2j-1)m +1}, \ldots,  \varepsilon_{2j m}) , j \geq 1 ) \, .
\eeq
Let $U_1 = \sum_{k=1}^m X_k$ and, for  any  integer $j  \in [2,N]$, define
\beq  \label{defUj}
U_j = \sum_{k=(2j-2) m+1}^{(2j-1)m}  ( X_{k,m} - \E ( X_{k,m}| {\mathbb F}_m )  )   \, . 
\eeq
For any integer $j \in [ 1,N]$, let 
\beq  \label{defRj}
R_j = \sum_{k=(2j-1)m+1 }^{2jm} ( X_{k,m}- \E ( X_{k,m}| {\mathbb F}_m )  )  \, ,
\eeq
\beq \label{notaSm1}
Y_j^{(1)}= U_j+ R_j \ \mbox{ and } \  S_{|m}^{(1)} = \sum_{j=1}^N Y_j^{(1)} \, .
\eeq
Let also 
\[
U_{N+1}  =  \sum_{k=2N m+1}^{ \min(n, (2N+1)m)}  ( X_{k,m} - \E (X_{k,m}| {\mathbb F}_m )  )
\]
and
\[R_{N+1}  = \sum_{k=(2N+1)m+1 }^{n} ( X_{k,m}- \E ( X_{k,m}| {\mathbb F}_m ))\, , \] where an empty sum has to be interpreted as $0$. 
Note that under ${\mathbb P}_{{\mathbb F}_m}$ (the conditional probability given  ${\mathbb F}_m$), the random vectors $(U_j,R_j)_{1 \leq j \leq N+1}$ are independent.  Moreover, by stationarity,  the r.v.'s 
$(U_j,R_j)_{2 \leq j \leq N}$ have the same distribution (as well as the   r.v.'s 
$(R_j)_{1 \leq j \leq N}$). 

Next, denoting by $ S_{|m}^{(2)} = \sum_{k=m+1}^n \E ( X_{k,m}| {\mathbb F}_m )$, the following decomposition is valid: 
\[
S_{n,m}:= \sum_{k=1}^m X_k +  \sum_{k=m+1}^{n} X_{k,m} =  S_{|m}^{(1)}  + S_{|m}^{(2)}  +  U_{N+1}  + R_{N+1} \, .
\]
{\it To simplify the exposition, assume in the rest of the proof that $n=2Nm$} (so that $m'=0$). There is no loss of generality by making such an assumption:  the only difference would be that  since $(U_{N+1}, R_{N+1}) $ does  not have the same law as the  $(U_j, R_j) $'s, $2 \leq j \leq N$, its  contribution  would have to be treated separately. 
 Therefore, from now we consider $m'=0$ and then  the following decomposition 
\beq \label{decSnm}
S_{n,m}= S_{|m}^{(1)}  + S_{|m}^{(2)}   \, .
\eeq
%
%\smallskip
%
%\noindent {\it Step 1. Outline of the proof.}  
We are now in position to give the main steps of the proof.  We start by writing 
\[
\big |  \E \big (  {\rm e}^{{\rm i} \xi  S_n /{\sqrt n}  }\big ) - {\rm e}^{-\xi^2s^2/2}  \big | 
\leq  \big |  \E \big (  {\rm e}^{{\rm i} \xi  S_{n}   /{\sqrt n}  }\big ) -\E \big (  {\rm e}^{{\rm i} \xi  S_{n,m}   /{\sqrt n}  }\big )  \big |  
 + \big |  \E \big (  {\rm e}^{{\rm i} \xi  S_{n,m} /{\sqrt n}  }\big ) - {\rm e}^{-\xi^2s^2/2}  \big |  \, .
\]
Next
\begin{multline*}
 \big |  \E \big (  {\rm e}^{{\rm i} \xi  S_{n,m} /{\sqrt n}  }\big ) - {\rm e}^{-\xi^2s^2/2}  \big |  \\
  =  \Big |  \E \Big (   {\rm e}^{{\rm i} \xi  S_{|m}^{(2)}  /{\sqrt n}  }  \Big [  \E_{{\mathbb F}_m} \big (  {\rm e}^{{\rm i} \xi  S_{|m}^{(1)}  /{\sqrt n}  } \big ) -  {\rm e}^{-\xi^2s^2/4}  \Big ] \Big ) 
+    {\rm e}^{-\xi^2s^2/4}  \Big ( \E \big (   {\rm e}^{{\rm i} \xi  S_{|m}^{(2)}  /{\sqrt n}  } \big )   - {\rm e}^{-\xi^2s^2/4} \Big )  \Big | \\
 \leq  \Big   \Vert  \E_{{\mathbb F}_m} \big (  {\rm e}^{{\rm i} \xi  S_{|m}^{(1)}  /{\sqrt n}  } \big ) -  {\rm e}^{-\xi^2s^2/4}  \Big  \Vert_{1}
+  \Big | \E \big (   {\rm e}^{{\rm i} \xi  S_{|m}^{(2)}  /{\sqrt n}  } \big )   - {\rm e}^{-\xi^2s^2/4} \Big |   \, .
\end{multline*}
Hence, starting from \eqref{IneBE} and selecting $T =1 /v_n  $ where $ \displaystyle v_n =   \Big (  \frac{\log n}{n }  \Big )^{q/2-1}$,  Inequality \eqref{ineBE1} of Theorem \ref{thmq=3} will follow if one can prove that 
\beq \label{inesmooth1}
 \int_{-T}^T   \frac{ \big |  \E \big (  {\rm e}^{{\rm i} \xi  S_{n}  /{\sqrt n}  }\big ) -\E \big (  {\rm e}^{{\rm i} \xi  S_{n,m}   /{\sqrt n}  }\big )  \big | }{|\xi|} {\rm d}\xi  \ll v_n \, ,
\eeq
\beq \label{inesmooth2}
 \int_{-T}^T   \frac{ \big   \Vert  \E_{{\mathbb F}_m} \big (  {\rm e}^{{\rm i} \xi  S_{|m}^{(1)}  /{\sqrt n}  } \big ) -  {\rm e}^{-\xi^2s^2/4}  \big  \Vert_{1} }{|\xi|} {\rm d}\xi  \ll v_n
\eeq
and
\beq \label{inesmooth3}
 \int_{-T}^T   \frac{ \big |  \E \big (   {\rm e}^{{\rm i} \xi  S_{|m}^{(2)}  /{\sqrt n}  } \big )   - {\rm e}^{-\xi^2s^2/4}   \big |  }{|\xi|} {\rm d}\xi  \ll v_n \, .
\eeq
The objective is then to  prove these three upper bounds, and the main differences compared to  \cite{Ji,Ji20} lie in the intermediate steps and the technical tools developed for this purpose. They will be based on the following dependence coefficients that are well adapted to our setting.  Let $p \geq 1$. For every $k\ge 1$, define
\beq\label{delta-def}
\delta_{p, \infty}^p (k) = \sup_{{\bar x}, {\bar y} \in X} \E \big |  X_{k,{\bar x} } - X_{k,{\bar y} }  \big |^p \, .
\eeq
If $\mu$ has a finite moment of order $q >1$, then,  by \cite[Prop. 3]{CDJ}, we know that   
\beq \label{estimatedelta}
\sum_{k \geq 1} k^{q-p-1} \, \delta_{p, \infty}^p (k) < \infty 
\qquad \forall p \in [1, q)\, .
\eeq
%and  
%\beq \label{estimatedelta2}
%\sum_{k \geq 1} k^{q-2} \, \delta_{p, \infty}^p (k) < \infty 
%\qquad \forall p \in (0,1] \, .
%\eeq
Hence, since $(\delta_{p, \infty} (k))_{k \geq 1}$ is 
non increasing, it follows that (if $\mu$ has a moment of order $q>1$) 
\beq \label{estimatedeltabis}\delta_{p, \infty} (k) = o \big ( 1/k^{q/p-1}  \big) 
\qquad \forall p\in [1,q)\, .
\eeq
%and 
%\beq \label{estimatedeltabis2}    \delta_{p, \infty} (k) = o \big (1/k^{(q-1)/p} \big ) 
%\qquad \forall p\in (0,1]\, .\eeq 
In the following commentary, we list the places where it is essential to use the coefficients $\delta_{p, \infty} (k)$ rather than the coupling coefficients used by  Jirak \cite{Ji20} in order to obtain the most accurate bounds possible. Note that this list is not exhaustive. 
\begin{Comment} \label{comment31}
Denote by $\vartheta'_k(p)$ and  $\vartheta^*_k(p)$ the coupling coefficients defined in \cite[Eq. (7)]{Ji20}. Note that in the Markovian case (which is our setting), these two coefficients are of the same order and can be bounded by 
$\delta_{p, \infty} (k) $. As  we shall see in  Lemma \ref{lmaR1normep}, by using a suitable Rosenthal-type inequality and the strength of the $\delta_{p, \infty}  $  coefficients, allowing to  control also the infinite norm of conditional expectations (see for instance \eqref{condexpectationrkm}),  we obtain, in particular,  $ \Vert R_1 \Vert_{p} \ll 1$  for   $p \geq 2 $ provided that $\mu$ has a moment of order $q=p+1$.  As a counterpart, Lemma 5.4 in \cite{Ji20} entails  that $\Vert R_1 \Vert_p \ll \sum_{k=1}^m  \delta_{p, \infty} (k) $,  and then $ \Vert R_1 \Vert_{p} \ll 1$  as soon as $\mu$ has a moment of order $q >2p$. A suitable control of   $ \Vert R_1 \Vert_{p} $  for some $p \geq 2$ is a key ingredient  to take care of the characteristic function of the $Y_j^{(1)} $'s conditionally to ${\mathbb F}_m$ that we will denote by  $\varphi_j(t)$ in what follows (see the definition \eqref{defvarphijx}). More precisely,  if  the condition (among others) $ \Vert R_1 \Vert_{p} \ll 1$ holds  for $p=2$, then we  get the upper  bound \eqref{conslma9} with $q=3$, and if it holds for $p=3$ then we get the better upper bound \eqref{q4varphij} (this difference in the upper bounds is the reason why in the statements of Theorem \ref{thmq=3} (with $q=3$) we have an extra logarithmic term compared to  Theorem \ref{thmq=4}).   Note that the upper bounds \eqref{conslma9} and \eqref{q4varphij}  come from Lemmas \ref{lma4.9},  \ref{lma4.9q=4} and \ref{lma4.9q=4bis}. Another crucial fact  that we would like to point out is the following: Imposing that $\mu$ has a moment of order $q=3$ implies   $ \Vert R_1 \Vert_{p} \ll 1$ only for $p=2$ and then, when $q \leq 3$,  Lemma 4.5 in \cite{Ji} cannot be used to get the upper bound \eqref{conslma4.5} which is essential to prove 
 \eqref{inesmooth2}. Indeed, in order for \cite[Lemma 4.5]{Ji} to be applied, it is necessary that $ \Vert R_1 \Vert_{p} \ll 1$ for  some $p>2$. The role of our Lemma \ref{lma4.5} is  then to overcome this drawback (see the step 2 below and in particular the control of both $ I_{1,N}(\xi)   $  and $ I_{3,N}(\xi) $). 
 
 On another hand, in view of \eqref{estimatedeltabis}, it is clear that, as $k \to \infty$,   for any $r \in [1, p[$, the coefficient $\delta_{r, \infty} (k) $ has a better behavior than $\delta_{p, \infty} (k)$. Hence, in some cases, it would be preferable to deal with the ${\mathbb L}^r$-norm rather than with the $ {\mathbb L}^p$-norm.  For instance, in our case, it is much more efficient to control $\Vert S_n  - S_{n,m} \Vert_1$ (see the forthcoming  upper bounds \eqref{step1P1S} and \eqref{step1P2S}) rather than $\Vert S_n  - S_{n,m} \Vert_p^p$ as done in Jirak \cite{Ji20} (see his upper bound (50)).   This is the reason why we can  start directly from Inequality \eqref{IneBE} and work with  the characteristic function rather  than using the decomposition given in \cite[Lemma 5.11]{Ji20}.
 %For instance if  the aim is to get in the Berry-Esseen theorem an upper bound of order $1/\sqrt{n} $ (up to some extra logarithmic terms) the condition $q >8$ required in \cite{Ji20} comes from this constraint, whereas we only need $q \geq 5/2$ at this stage.  This is the reason why we can  start directly from Inequality \eqref{IneBE} and work with  the characteristic function rather  than using the decomposition given in \cite[Lemma 5.11]{Ji20}.
% In both cases these quantities have to be controlled by $1/\sqrt{n}$  and to see the differences between the two upper bounds take $m$ equals to $n$ up to a logarithmic term both in \eqref{step1P2S} and in \cite[Ineq. (50)]{Ji20}.  
 \end{Comment}
 
 \smallskip
 
Let us now come back to the proof. The next steps will consist in proving the upper bounds \eqref{inesmooth1}-\eqref{inesmooth3}. 

\medskip

\noindent \textit{Step 1. Proof of \eqref{inesmooth1}.}  Note that 
\[
 \int_{-T}^T   \frac{ \big |  \E \big (  {\rm e}^{{\rm i} \xi  S_n  /{\sqrt n}  }\big ) -\E \big (  {\rm e}^{{\rm i} \xi  S_{n,m}   /{\sqrt n}  }\big )  \big | }{|\xi|} {\rm d}\xi  \leq \frac{T}{\sqrt{n}}\Vert S_n  - S_{n,m} \Vert_{1}  \, .
\]
But,  by stationarity and  \cite[Lemma 24]{CDM} (applied with  $M_k = +\infty$), 
\begin{equation} \label{step1P1S}
\Vert S_n  - S_{n,m} \Vert_1 \leq n \Vert X_{m+1}-X_{m+1,m} \Vert_1 \leq n \delta_{1, \infty} (m) \, .
\end{equation}
Hence, by \eqref{estimatedeltabis} and the fact that $\mu$ has a moment of order $q>1$, we derive
\begin{equation}\label{step1P2S}
\Vert S_n - S_{n,m} \Vert_1  \ll n m^{-(q-1)}  \, .
\end{equation}
So, overall, since  $T \ll m^{q/2-1}$, it follows that 
\[
 \int_{-T}^T   \frac{ \big |  \E \big (  {\rm e}^{{\rm i} \xi  S_n  /{\sqrt n}  }\big ) -\E \big (  {\rm e}^{{\rm i} \xi  S_{n,m}   /{\sqrt n}  }\big )  \big | }{|\xi|} {\rm d}\xi  \ll  \frac{n^{1/2}}{m^{q/2}}  \, .
\]
The upper bound \eqref{inesmooth1} follows from the fact that  $m \sim \kappa_2 n (\log  n)^{-1}$.

\medskip

\noindent \textit{Step 2. Proof of \eqref{inesmooth2}.} For  any $x\in {\mathbb R}$ and any integer $j \in [1,N]$, let 
\beq \label{defvarphijx}
\varphi_j (x)= \E \Big ( {\rm e}^{{\rm i} x Y_j^{(1)} / \sqrt{ 2  m} }  | {\mathbb F}_m\Big )   \, . \eeq
Since, under ${\mathbb P}_{{\mathbb F}_m}$, the $Y_j^{(1)} $'s are independent we write
\beq \label{ineBEclassic}
 \big   \Vert  \E_{{\mathbb F}_m} \big (  {\rm e}^{{\rm i} \xi  S_{|m}^{(1)}  /{\sqrt n}  } \big ) -  {\rm e}^{-\xi^2s^2/4}  \big  \Vert_{1}  =  \E \Big [   \Big | \prod_{j=1}^N  \varphi_j \Big (\frac{ \xi}{\sqrt{N}} \Big )   - \prod_{j=1}^N  {\rm e}^{- \xi^2s^2/(4N) } \Big | \Big ]  \, .
\eeq
 As in \cite[Section 4.1.1]{Ji}, we use the following basic identity: for 
 any complex numbers $(a_j)_{1\le j\le N}$ and $(b_j)_{1\le j\le N}$,
$
  \prod_{j=1}^N  a_j  -  \prod_{j=1}^N b_j = \sum_{i=1}^n ( \prod_{j=1}^{i-1} b_j)  (a_i-b_i)  ( \prod_{j=i+1}^{N} a_j) 
$ to handle the right-hand side of \eqref{ineBEclassic}. Taking into account that $(\varphi_j(t))_{1 \leq j \leq N}$ forms a one-dependent sequence and that  the r.v.'s 
$(U_j,R_j)_{2 \leq j \leq N}$ have the same distribution, we then infer that  
\begin{equation} \label{ineBE4P0}
\E \Big [   \Big | \prod_{j=1}^N  \varphi_j \Big (\frac{ \xi}{\sqrt{N}} \Big )   -  \prod_{j=1}^N  {\rm e}^{- \xi^2/(4N) }\Big | \Big ]  \leq I_{1,N}(\xi) +I_{2,N}(\xi) + I_{3,N}(\xi)   \, ,
\end{equation}
where 
\[
 I_{1,N}(\xi)   = (N  -1) \Vert \varphi_2 (\xi /{\sqrt N})  -  {\rm e}^{- \xi^2s^2/(4N) } \Vert_{1} \Big \Vert    \prod_{j=N/2}^{N-1}  \Big |  \varphi_j \Big (\frac{ \xi}{\sqrt{N}} \Big )  \Big |   \Big \Vert_{1} \, , 
\]
\[
 I_{2,N}(\xi) = N   {\rm e}^{- \xi^2s^2 (  N-6) /(8N) }   \Vert \varphi_2 (\xi /{\sqrt N})  -  {\rm e}^{- \xi^2s^2/(4N) } \Vert_{1}  \]
and
\[
 I_{3,N}(\xi) =     \Vert \varphi_1 (\xi /{\sqrt N})  -  {\rm e}^{- \xi^2s^2/(4N) } \Vert_{1} \Big \Vert    \prod_{j=N/2}^{N-1}  \Big |  \varphi_j \Big (\frac{ \xi}{\sqrt{N}} \Big )  \Big |   \Big \Vert_{1}  \, . \]
To integrate the above quantities,  we need to give suitable upper bounds for the two terms $\Vert \varphi_j (t)  -  {\rm e}^{-s^2 t^2/ 4 } \Vert_{1}$ and $\Vert \prod_{j=N/2}^{N-1} |  \varphi_j (t) | \Vert_{1}$. Applying the first part of Lemma \ref{lma4.9} and using stationarity, we derive that for any $2 \leq j \leq N$, 
\beq \label{conslma9}
\Vert \varphi_j (t)  -  {\rm e}^{- s^2 t^2/ 4 } \Vert_{1} \ll  \frac{ t^2 }{m^{q/2 -1}}  + \frac{|t|}{m^{q-3/2}} \, .
\eeq
Moreover the second part of Lemma \ref{lma4.9} implies that 
\beq \label{conslma9j=1}
\Vert \varphi_1 (t)  -  {\rm e}^{- s^2 t^2/ 4 } \Vert_{1} \ll \frac{ t^2 }{m^{q/2 -1}}   \, .
\eeq

\noindent On another hand, according to \cite[Inequality (4.14)]{Ji}, for any integer $\ell \in [1, m]$, 
\[
\Big \Vert \prod_{j=N/2}^{N-1} |  \varphi_j (t) | \Big  \Vert_{1} \leq \Big  \Vert \prod_{j \in {\mathcal J}}  \big |  \varphi^{(\ell)}_j (t\sqrt{(m-\ell )/(2m)}) \big  | \Big  \Vert_{1}  \, , 
\]
where ${\mathcal J} = [N/2, N-1] \cap 2 {\mathbb N}$, 
\[
 \varphi^{(\ell)}_j ( x ) = \E\Big ( {\rm e}^{{\rm i} x H^{(\ell)}_{j,m}   }  \big | {{\mathcal H}^{(\ell)}_{j,m}}  \Big ) 
\]
with ${\mathcal H}^{(\ell)}_{j,m}= {\mathbb F}_m \vee  \sigma( \varepsilon_{ 2(j-1)m +1}, \ldots,  \varepsilon_{ 2(j-1)m +\ell}) $ and 
\[
H^{(\ell)}_{j,m}  =  \frac{1}{\sqrt{m- \ell}} \Big (  \sum_{k=2(j-1)m + \ell +1}^{(2j-1)m}( X_{k,m} - \E ( X_{k,m}| {\mathcal H}^{(\ell)}_{j,m} )  )   + R_j - \E (R_j| {\mathcal H}^{(\ell)}_{j,m} )   \Big ) \, .
\]
We shall apply Lemma \ref{lma4.5} with 
\[A_j = \frac{1}{\sqrt{m- \ell}}  \sum_{k=2(j-1)m + \ell +1}^{(2j-1)m}( X_{k,m} - \E ( X_{k,m}| {\mathcal H}^{(\ell)}_{j,m} )  ), \ B_j = \frac{R_j -  \E (R_j| {\mathcal H}^{(\ell)}_{j,m} ) }{m^{(3-q)/2}}  \]
and 
$\displaystyle  a = \frac{ m^{(3-q)/2} }{(m- \ell)^{1/2} } $. By stationarity, for any $j \in {\mathcal J}$, 
\begin{multline*}
{\mathbb P}  \big ( \E_{H^{(\ell)}_{j,m} } (A^2_j  )  \leq s^2/4 \big  ) = {\mathbb P} \big ( \E_{H^{(\ell)}_{2,m} } (A^2_2  )  \leq s^2/4 \big )  \\
=  {\mathbb P} \Big ( (m-\ell)^{-1} \E_{m} \Big  ( \Big ( \sum_{k=m+1}^{2m-\ell}  (X_{k,m} - \E_m ( X_{k,m} ) \Big )^2 \Big )  \leq s^2/4  \Big ) \, ,
\end{multline*}
where $\E_m (\cdot)$ means $\E(\cdot |  {\mathcal G}_m )$ with ${\mathcal G}_m = \sigma (W_0, \varepsilon_1, \ldots, \varepsilon_m)$.  Let $K$ be a positive integer and note that
\begin{align*}
& \Big |  \Big \Vert  \sum_{k=m+1}^{m+K }  (X_{k,m} - \E_m ( X_{k,m} ) ) \Big \Vert_2 - \Big  \Vert  \sum_{k=m+1}^{m+K }  X_{k} \Big   \Vert_2 \Big |   \\ &  \quad \quad  \leq 
 \Big \Vert  \sum_{k=m+1}^{m+K }  ( X_{k,m} - X_k )  \Big \Vert_2 + \sum_{k=m+1}^{m+K }  \Vert  \E_m ( X_{k,m} ) \Vert_{\infty} \, .
\end{align*}
But, by using the remark after \cite[Prop. 3]{CDJ}, we infer that, for $k \geq m+1$, 
\beq \label{Borne1condexpect}
  \Vert  \E_m ( X_{k,m} ) \Vert_{\infty}   \leq  \delta_{1, \infty} (k-m)  \, .
\eeq
Next,  by \cite[Lemma 24]{CDM} (applied with  $M_k = +\infty$), 
\begin{align*}
\Big \Vert  \sum_{k=m+1}^{m+K }  ( X_{k,m} & - X_k )  \Big \Vert^2_2  =  \sum_{k=m+1}^{m+K }  \Vert X_{k,m} - X_k \Vert_2^2 \\ &  \quad \quad + 2  \sum_{k=m+1}^{m+K-1 }   \sum_{\ell=k+1}^{m+K } 
\E \Big (  (X_{k,m} - X_k) \E_{k}  (X_{\ell,m} - X_\ell)   \Big )   \\
& \leq    K \delta^2_{2 , \infty} (m)+  2  \sum_{k=m+1}^{m+K-1 }   \sum_{\ell=k+1}^{m+K }   \Vert \E_{k}  (X_{\ell,m} - X_\ell)  \Vert_{\infty}   \Vert X_{k,m} - X_k \Vert_1 \\
& \leq    K \delta^2_{2 , \infty} (m)+  2  \sum_{k=m+1}^{m+K-1 }   \sum_{\ell=k+1}^{m+K }  \delta_{1, \infty} (\ell - k)   \delta_{1, \infty}  (m) \, .
\end{align*}
Therefore, by taking into account \eqref{estimatedeltabis} and the fact that $\mu$ has a moment of order $q  >2$, we get that 
\[
\Big \Vert  \sum_{k=m+1}^{m+K }  ( X_{k,m}  - X_k )  \Big \Vert^2_2  = o (  K m^{2-q } ) \, ,
\]
which combined with \eqref{Borne1condexpect} implies that
\[
K^{-1/2}\Big | \Big \Vert  \sum_{k=m+1}^{m+K }  (X_{k,m} - \E_m ( X_{k,m} ) ) \Big \Vert_2 - \Big  \Vert  \sum_{k=m+1}^{m+K }  X_{k} \Big   \Vert_2 \Big |   \ll  m^{1-q/2 }  + K^{-1/2} \, .
\]
But, using stationarity,  $K^{-1/2}  \big  \Vert  \sum_{k=m+1}^{m+K }  X_{k} \big   \Vert_2 = K^{-1/2}  \big  \Vert  \sum_{k=1}^{K }  X_{k} \big   \Vert_2  \rightarrow s >0 $.  Hence provided that $(m-\ell) $ is large enough, we have 
\beq  \label{restrictioonell}
 (m-\ell)^{-1} \E \Big  ( \Big ( \sum_{k=m+1}^{2m-\ell}  (X_{k,m} - \E_m ( X_{k,m} )  )\Big )^2 \Big )  > s^2/2\, . 
\eeq
So, overall,  setting ${\bar X}_{k,m} := X_{k,m} - \E_m ( X_{k,m} ) $,  for $(m-\ell) $ large enough,  we get 
\begin{multline*}
{\mathbb P} \big ( \E_{H^{(\ell)}_{2,m} } (A^2_2  )  \leq s^2/4 \big )  \\
\leq   {\mathbb P} \Big ( (m- \ell)^{-1} \Big | \E_{m} \Big  ( \Big ( \sum_{k=m+1}^{2m-\ell} {\bar X}_{k,m} \Big )^2  -  \E \Big  ( \Big ( \sum_{k=m+1}^{2m-\ell}  {\bar X}_{k,m}  \Big )^2\Big )  \Big |  \geq \frac{s^2}{4 }  \Big ) \, .
\end{multline*}
Using Markov's inequality,   the same arguments as those used in the proof of Lemma \ref{lma4.7}, and since $q >2$, we then derive that, for $(m-\ell) $ large enough and  any $j \in {\mathcal J}$, 
\[
{\mathbb P}  \big ( \E_{H^{(\ell)}_{j,m} } (A^2_j  )  \leq s^2/4 \big  )  \ll (m - \ell)^{-\varepsilon} \text{ for some $\varepsilon >0$.} 
\]
Hence, provided that $m-\ell$ is large enough, Item (ii) of Lemma \ref{lma4.5} is satisfied with $u^- = s^2/4$.  Note now that by stationarity, for any $j \in {\mathcal J}$, 
\[
\E(B_j^2)  \leq 4  \frac{\E (R_j^2)  } {m^{3-q}}= 4 \frac{\E (R_1^2) } {m^{3-q}} \ll 1 \, , 
\]
by using Lemma \ref{lmaR1normep} with $p=2$. This proves Item (iv) of Lemma \ref{lma4.5}. Next, for $p\geq 2$, using stationarity and \cite[Cor. 3.7]{MPU19}, we get that for any $j \in {\mathcal J}$, 
\beq \label{momentAJ}
\E(|A_j|^p)  \leq 2^p (m- \ell)^{-p/2}  \Big \Vert \sum_{k=m+1}^{2m - \ell} X_{k,m}  \Big  \Vert_{p} \ll \Big [  \Vert X_{1+m,m} \Vert_{p} + \sum_{k=m+1}^{2m-\ell} k^{-1/2} \Vert \E_m (X_{k,m}) \Vert_{p}  \Big ]^p  \, .
\eeq
But $\Vert X_{1+m,m} \Vert_{p}  \leq \Vert X_{1} \Vert_{p}  < \infty$ if $p \leq q$ (indeed recall that it is assumed that $\mu$ has a moment of order $q$) and, by \eqref{Borne1condexpect},  $\Vert \E_m (X_{k+m,m}) \Vert_{p} \leq \delta_{1, \infty} (k)$. Hence, by \eqref{estimatedelta}  and since $\mu$ has a moment of order $q>2$, Item (iii) of Lemma \ref{lma4.5} is satisfied for $p=q$.  So, overall, noticing that $|{\mathcal J}| \geq N/8\ge 16$,  we can apply Lemma \ref{lma4.5} to derive that there exist positive finite constants $c_1$, $c_2$ and $c_3$ depending in particular on $s^2$ but not on $(m,n)$ such that for $(m-\ell)$ large enough (at least such that $a = \frac{ m^{(3-q)/2} }{(m- \ell)^{1/2} }  \leq c_1$), we have
\[
 \Big  \Vert \prod_{j \in {\mathcal J}}  \big |  \varphi^{(\ell)}_j (x) \big  | \Big  \Vert_{1}  \le {\rm e}^{-c_3 x^2 N/8} + {\rm e}^{-  N/256}\  \   \text{for $x^2 \leq c_2$, }  
\]
implying overall that, for $(m-\ell)$ large enough and for $t^2 (m-\ell)/(2m) \leq c_2$,    
\beq \label{conslma4.5}
\Big \Vert \prod_{j=N/2}^{N-1} |  \varphi_j (t) | \Big  \Vert_{1}  \le {\rm e}^{-c_3 t^2 (m-\ell) N/ ( 16 m) } + {\rm e}^{- N/256}  \, . 
\eeq
The  bounds \eqref{conslma9},  \eqref{conslma9j=1} and \eqref{conslma4.5} allow to give an upper bound for  the terms $I_{1,N}(\xi)$, $I_{2,N}(\xi)$ and $I_{3,N}(\xi)$ and next to integrate them over $[-T,T]$ when they are divided by $|\xi |$. Hence the computations in \cite[Sect. 4.1.1., Step 4]{Ji} are replaced by the following ones. First, as in \cite{Ji}, we select 
\beq \label{selectl}
\ell = \ell(\xi) =  {\bf 1}_{ \{\xi^2 < N c_2\}} + ( m - [nc_2/(2\xi^2)] +1 ) {\bf 1}_{\{\xi^2 \geq  N c_2\}}  \, . 
\eeq
Therefore $m - \ell $ is either equal to $m-1$ or to $  [nc_2/(2\xi^2)] -1 $. Since $|\xi | \leq T = \big ( n / (\log n ) \big )^{q/2-1}$, it follows that  $nc_2/(2\xi^2) \geq 2^{-1}c_2  (\log n)^{q-2} n^{3-q }$.Therefore
\[
a = \frac{ m^{(3-q)/2} }{(m- \ell)^{1/2} }  \leq   \frac{ m^{(3-q)/2} }{m-1 }    +     \frac{  2 m^{(3-q)/2} }{ c_2 n^{3-q } ( \log n )^{q-2} }  \, , 
\]
which is going to zero as $n \rightarrow \infty$ by the selection of $m$. Then, for any $c_1 >0$, we have $a <c_1$ for $n $ large enough.  This justifies the application of Lemma \ref{lma4.5}. 
So,  starting from \eqref{conslma4.5} and taking into account the selection of $\ell$, we get that for any $|\xi| \leq T$ and $n$ large enough, 
\beq \label{conslma4.5bis}
\Big \Vert \prod_{j=N/2}^{N-1} |  \varphi_j (\xi/{\sqrt N}) | \Big  \Vert_{1}  \ll {\rm e}^{-c_3 \xi^2 /32} {\bf 1}_{ \{\xi^2 < N c_2\}}  + {\rm e}^{- c_3 c_2 N/32 } {\bf 1}_{ \{\xi^2 \geq N c_2\}}  + {\rm e}^{- N/256}  \, . 
\eeq
Select now 
\beq  \label{selectionofN}
N  = [\kappa \log n  ] \ \  \text{with $\kappa >2 \max( 256, 32 (c_2c_3)^{-1})  $ } 
\eeq 
and then $m \sim  (2 \kappa )^{-1}n / \log n$. Taking into account \eqref{conslma9}, \eqref{conslma9j=1} and \eqref{conslma4.5bis},   we get, for $n $ large enough, 
\begin{multline} \label{ineBE4P1}
\int_{-T}^T  ( I_{1,N}(\xi) +  I_{3,N}(\xi)  ) / |\xi | \, {\rm d}\xi \ll  N  \int_{0}^T \Big ( \frac{|\xi|}{ N m^{q/2-1}} + \frac{1}{ \sqrt{N } m^{q-3/2}}  \Big )  \Big (  {\rm e}^{-c_1 \xi^2 /32}  +n^{-2} \Big )  \, {\rm d}\xi  \\
\ll   \frac{1}{m^{q/2-1} } +   \frac{\sqrt{N}}{ m^{(q-1)/2}m^{q/2-1}}   +  \frac{ T^2 }{n^2 m^{q/2-1}} +   \frac{ T \sqrt{N}}{ n^2 m^{(q-1)/2}m^{q/2-1} }  \ll  \Big ( \frac{\log n }{n} \Big )^{q/2-1} \, .
\end{multline}
Next,  using \eqref{conslma9}, we derive 
\begin{align*} 
I_{2,N}(\xi)   \ll   \Big ( \frac{\xi^2}{m^{q/2-1} }  + \frac{ \sqrt{N}|\xi |}{m^{q-3/2}}   \Big )  \times \  {\rm e}^{- s^2 \xi^2/  16 }   \, .
\end{align*}
Therefore, by the selection of $m$ and $N$, 
\begin{equation}  \label{ineBE4P2}
\int_{-T}^T I_{2,N}(\xi) /  |\xi | \, {\rm d}\xi  \ll   \Big ( \frac{\log n }{n} \Big )^{q/2-1}    \, .
\end{equation}
Starting from \eqref{ineBEclassic} and taking into account \eqref{ineBE4P0}, \eqref{ineBE4P1} and  \eqref{ineBE4P2}, the upper bound in \eqref{inesmooth2}  follows. 

\medskip

\noindent \textit{Step 3. Proof of \eqref{inesmooth3}.} Recall that $S_{|m}^{(2)} = \sum_{k=m+1}^n \E ( X_{k,m}| {\mathbb F}_m )$, and  recall  that we assume that $2Nm = n$. Denoting 
\[  Y_j^{(2)}  =  U_j^{(2)}  +  R_j^{(2)} \ \ \text{for} \ j =1, \ldots, N \, ,  \] 
where $ U_N^{(2)} =  \sum_{k=(2N-1)m+1}^{n} \E ( X_{k,m}| {\mathbb F}_m )$, $R_N^{(2)} =0$, 
\[
 U_j^{(2)}  =  \sum_{k=(2j-1)m+1}^{2jm} \E ( X_{k,m}| {\mathbb F}_m ) \,  \mbox{ and }  \, R_j^{(2)}  =  \sum_{k=2jm+1}^{(2j+1)m} \E ( X_{k,m}| {\mathbb F}_m )\ \  \text{for}\  j =1, \ldots, N-1 \, , 
\]
we have $S_{|m}^{(2)} = \sum_{j=1}^{N} Y_j^{(2)}$.  Note that the random vectors $( U_j^{(2)} ,R_j^{(2)} )_{1 \leq  j \leq N}$ are independent.  The proof of   \eqref{inesmooth3} can be done by using similar (but even simpler) arguments to those developed in the step 2. In this part, one of the important fact is to notice that the $R^{(2)}_j$'s also have a negligible contribution. Indeed,  for any 
$2m+1 \leq k \leq 3m$, 
\begin{multline*}
\Vert \E ( X_{k,m}| {\mathbb F}_m ) \Vert_{\infty} 
= \Big  \Vert  \int \! \! \int \Big ( f_m ( \varepsilon_{k-m+1} , \ldots, \varepsilon_{2m}, a_{2m+1}, \ldots, a_{k})  \\ - f_m ( b_{k-m+1} , \ldots, b_{2m}, b_{2m+1}, \ldots, b_{k}) \Big )  \prod_{i=2m+1}^{k} d \mu (a_i)  \prod_{i=k-m+1}^{k} d \mu (b_i) \Big  \Vert_{\infty} \\
\leq \sup_{{\bar x}} \Big | \E (X_{k-2m} | W_0={\bar x}) - \int \E (X_{k-2m} | W_0={\bar y}) d\nu ({\bar y})  \Big | \leq \delta_{1, \infty} (k-2m) \, .
\end{multline*}
Hence by stationarity,  \eqref{estimatedelta} and since $q \geq 2$, we derive that $\Vert R^{(2)}_j \Vert_{ \infty} \ll 1$ for any $j=1, \ldots, N$.

\medskip

To complete the proof of the upper bound \eqref{ineBE1}, we just have to put together the results in the steps 1, 2 and 3.  \qed

\subsubsection{Proof of the upper bound  \eqref{ineBE1bis}} \label{subsection2} Recall the notation $S_{n,{\bar u}}:= \sum_{k=1}^n X_{k,{\bar u}}$ where $X_{k,{\bar u}}$ denotes the random variable $X_k$  defined by \eqref{defMCXk} when the Markov chain $(W_n)_{n \geq 0}$ starts from $\bar u \in X$.   Our starting point is the following upper bound: 
\begin{equation}\label{lemme-bougerol}
\sup_{n \geq 1} \Big \Vert  \log (\Vert A_n \Vert) -n \lambda_\mu -  \int_X S_{n,{\bar u}} d  \nu ({\bar u})   \Big \Vert_{\infty} < \infty \, .
\end{equation}
The proof of \eqref{lemme-bougerol} is outlined in  Section 8.1 in \cite{CDJ} but,  since it is a key ingredient in the proof of \eqref{ineBE1bis}, we shall provide more details here.    
Let $g\in G$ and $\bar u\in X$. By item $(i)$ of 
Lemma 4.7 in \cite{BQ}, there exists $\bar v(g)$ such that 
\[
\log \|g\|- \sigma(g,\bar u)\le -\log  \delta \big(\bar u, \bar v(g)\big)   \, ,
\]
where $\delta (\bar u,\bar v):= \frac{|\langle u,v\rangle|}{\|u\|\, \|v\|}$. 
Integrating with respect to $\nu$, it follows that 
\begin{equation}\label{BQ-est}
0\le \log \|g\| -\int_X \sigma(g, \bar u)\, d\nu(\bar u)\le  \sup_{\bar v\in X} \int_X|\log \delta(\bar u, \bar v)|\, d\nu(\bar u) \, .
\end{equation}
But, according to  Proposition 4.5 in \cite{BQ}, since $\mu$ has a polynomial moment of order $q \geq 2$, $ \sup_{\bar v \in X}\int_X \big| \log \delta({\bar u}, \bar v)\, \big|\, d  \nu ({\bar u}) < \infty$.  Therefore, 
\eqref{lemme-bougerol} follows from an application of  \eqref{BQ-est} with $g= A_n$.

Now, using \eqref{lemme-bougerol} and Lemma 1 in \cite{Bo82},  the upper bound  \eqref{ineBE1bis} will follow if one can prove that 
\beq \label{ineBE1bispr}
  \sup_{y \in {\mathbb R}} \Big | {\mathbb P} \Big ( \int_X S_{n,{\bar u}} d  \nu ({\bar u})    \leq y \sqrt{n} \Big ) - \Phi (y/ s)  \Big | \ll \Big ( \frac{\log n }{n} \Big )^{q/2-1}   \, .
\eeq
We proceed as in the proof of the upper bound  \eqref{ineBE1} with the following differences. First we consider 
\[
S_{n,m} = \sum_{k=1}^m \int_X  X_{k,{\bar u}} d  \nu ({\bar u})  + \sum_{k=m+1}^n X_{k,m} \, ,
\]
where $X_{k,m}$ is defined by \eqref{defXkm}. 
Hence 
\[
\Big \Vert  \int_X S_{n,{\bar u}} d  \nu ({\bar u})  - S_{n,m} \Big  \Vert_1 \leq  \int_X \sum_{k=m+1}^n \Vert   X_{k,{\bar u}}-  X_{k,m} \Vert_1 d  \nu ({\bar u}) \leq n \delta_{1,\infty} (m)  \, .
\]
It follows that the step 1 of the previous subsection is unchanged. Next, we use the same notation as in Subsection \ref{subsection1} with the following change: $U_1$ is now defined by
\beq \label{defU1bis}
U_1 = \sum_{k=1}^m \int_X  X_{k,{\bar u}} d  \nu ({\bar u})  \, ,
\eeq
and then, when $n =2mN$,  the decomposition \eqref{decSnm} is still valid for 
$S_{n,m}$.  The step 3 is also unchanged. Concerning the step 2, the only difference concerns  the upper bound of the quantity $\Vert \varphi_1 (t)  -  {\rm e}^{- s^2 t^2/ 4 } \Vert_{1} $ since the definition of $U_1$ is now given by \eqref{defU1bis}.  To handle this term, we note that for $f(x) \in \{ \cos x, \sin x \}$,  by using the arguments used in the proof of \cite[Lemma 24]{CDM}, 
we have 
\begin{multline*}
\Big \Vert \E_{{\mathbb F}_m}  \Big [ f \Big ( t \frac{  \sum_{k=1}^m \int_X  X_{k,{\bar u}} d  \nu ({\bar u})  +R_1 }{\sqrt{2m}} \Big ) \Big ]  -  \E_{{\mathbb F}_m}  \Big [ f \Big ( t \frac{  \sum_{k=1}^m  X_{k}   +R_1 }{\sqrt{2m}} \Big ) \Big ]   \Big   \Vert_{1} 
\\ \leq  \frac{|t|}{\sqrt{2m}}  \int_X \sum_{k=1}^m \Vert   X_{k,{\bar u}}-  X_{k} \Vert_1 d  \nu ({\bar u}) \leq  \frac{|t|}{\sqrt{2m}} \sum_{k= 1}^m \delta_{1, \infty} (k)  \ll \frac{|t|}{\sqrt{m}}  \, .
\end{multline*}
The last upper bound follows from \eqref{estimatedelta} together with the fact that $\mu$ is assumed to have a moment of order at least $2$. 
Next, by taking into account \eqref{conslma4.5bis}, note that 
\[
\int_{-T}^{T}  \frac{|\xi|}{ \sqrt{N}\sqrt{m}}  \Big \Vert \prod_{j=N/2}^{N-1} |  \varphi_j (\xi/{\sqrt N}) | \Big  \Vert_{1} d \xi  \ll 1/\sqrt{n} \, .
\]
This implies in particular that  \eqref{ineBE4P1} still holds.   Compared to Subsection  \ref{subsection1}  the rest of the proof is unchanged. \qed

\subsubsection{Proof of the upper bound  \eqref{ineBE1ter}}  \label{subsection3} 

Once again we highlight the differences with respect to the proof given in Subsection \ref{subsection1}. For $x \in S^{d-1}$, we consider 
\[
S_{n,m,\bar x} = \sum_{k=1}^m X_{k,\bar x}  + \sum_{k=m+1}^n X_{k,m} \, ,
\]
and we note that \[
\sup_{\bar x \in X} \Vert S_{n,\bar x} -S_{n,m,\bar x} \Vert_1 \leq \sum_{k=m+1}^n \sup_{\bar x\in X} \Vert  X_{k,\bar x}-  X_{k,m} \Vert_1  \leq n \delta_{1,\infty} (m)  \, .
\]
Once again  Step 1 of Subsection \ref{subsection1} is unchanged. Next, $U_1$ is now defined by
\beq \label{defU1ter}
U_{1,\bar x}  =U_1= \sum_{k=1}^m  X_{k,\bar x}  \, ,
\eeq
and the step 3 is also unchanged. Concerning the step 2,  due to the new definition \eqref{defU1ter} of $U_1$, the only difference concerns again the upper bound of the quantity $\Vert \varphi_{1} (t)  -  {\rm e}^{- s^2 t^2/ 4 } \Vert_{1} $.  To handle this term, we note that for $f(y) \in \{ \cos y, \sin y \}$,  
we have, by using \eqref{estimatedelta} together with the fact that $\mu$ is assumed to have a moment of order at least $2$, 
\begin{multline*}
\sup_{\bar x\in X}\Big \Vert \E_{{\mathbb F}_m}  \Big [ f \Big ( t \frac{  \sum_{k=1}^m \ X_{k,\bar x}  +R_1 }{\sqrt{2m}} \Big ) \Big ]  -  \E_{{\mathbb F}_m}  \Big [ f \Big ( t \frac{  \sum_{k=1}^m  X_{k}   +R_1 }{\sqrt{2m}} \Big ) \Big ]   \Big   \Vert_{1} 
\\ \leq  \frac{|t|}{\sqrt{2m}}  \sum_{k=1}^m  \sup_{\bar x  \in  X}\Vert  X_{k,\bar x}-  X_{k} \Vert_1  \leq  \frac{|t|}{\sqrt{2m}} \sum_{k= 1}^m \delta_{1, \infty} (k)  \ll \frac{|t|}{\sqrt{m}} \, .
\end{multline*}
We then end the proof as in Subsection \ref{subsection2}.  \qed

%\[
%\Big \Vert \E_{{\mathbb F}_m}  \Big [ f \Big ( t \frac{Y_1^{(1)}}{\sqrt{2m}} \Big ) \Big ]  -  \E \big [ f  ( t  s N /{\sqrt 2} ) \big ] \Big   \Vert_{1}  \ll |t|^3 m^{-1/2} + |t| m^{-1} (\log m) \, .
%\]

\subsection{ Proof of Theorem \ref{thmq=4}}

Let us point out the  differences compared to the proof of Theorem \ref{thmq=3} (the selections of $N$ and $m$ being  identical). To get the upper bound \eqref{conslma4.5bis}, we still establish an upper bound similar to 
\eqref{conslma4.5}  valid for any $\ell \in [1, m]$ and any $t$ such that $t^2 (m-\ell)/ (2m) \leq C$ for some positive constant $C$. Since $\mu$ has a finite moment of order $q=4$, according to Lemma \ref{lmaR1normep},  $ \Vert R_1 \Vert_{3} \ll 1$. Hence, using  Lemma 
\ref{lma4.5}  with 
\[A_j = \frac{1}{\sqrt{m- \ell}}  \Big (  \sum_{k=2(j-1)m + \ell +1}^{(2j-1)m}( X_{k,m} - \E ( X_{k,m}| {\mathcal H}^{(\ell)}_{j,m} )  ) + R_j -  \E (R_j| {\mathcal H}^{(\ell)}_{j,m} )  \Big )  \]
and $a=0$ (here  Lemma 4.5 in \cite{Ji} can also be used), the desired upper bound follows and the constant $C$ appearing above in the restriction for $t$ can be taken equal to $c_2$ (which is the constant appearing in Lemma 
\ref{lma4.5}). The fact that $a=0$ implies that we do not need to verify, as in the proof of Theorem \ref{thmq=3}, that $m^{(3-q)/2} (m- \ell)^{-1/2} \leq c_1$.  Next, we select $\ell$ as in \eqref{selectl}. This selection makes sense if $\xi^2 \leq n c_2 /2$. Therefore, we use \eqref{IneBE} by  selecting
$
T = \eta \sqrt{n}$ with $\eta$ small enough (more precisely such that $ c_2 /(2 \eta^2) $ is large enough for \eqref{restrictioonell} to be satisfied when $m-\ell$ is of order $ c_2 /(2 \eta^2) $). Therefore,  for any $|\xi| \leq T$, the upper bound \eqref{conslma4.5bis} is still valid.  The second difference, in addition to the choice of $T$, is that  instead of using Lemma \ref{lma4.9},  we use Lemmas \ref{lma4.9q=4} and \ref{lma4.9q=4bis} with $r=3$  which then entail that for any $j \geq 1$, 
\beq \label{q4varphij}
\Vert \varphi_j (\xi / \sqrt{N})  -  {\rm e}^{- s^2 \xi^2/ (4N) } \Vert_{1} \ll  N^{-1} |\xi|^3 n^{-1/2}  + |\xi | n^{-1/2} m^{-3/10}  \, .  \quad \qed
\eeq
Note that the upper bound (42)  in Jirak \cite{Ji20} with $p=3$ has the same order as \eqref{q4varphij} and is obtained provided $\sum_{k \geq 1} k^a \delta_{3, \infty } (k) < \infty$ for some $a>0$ (indeed \cite[Lemma 5.8 (iii)]{Ji20} is a key ingredient to get (42)). Now, using \eqref{estimatedelta}, we see that $\sum_{k \geq 1} k^a \delta_{3, \infty } (k) < \infty$ for some $a>0$   as soon as $\mu$ has a moment of order $q >6$. Actually \cite[Lemma 5.8]{Ji20} is not needed in its full generality to get an upper bound as  \eqref{q4varphij}. Indeed  our Lemmas \ref{lma4.9q=4} and \ref{lma4.9q=4bis} are rather based on an estimate as  \eqref{zolotarev3-step3P2} which involves the ${\mathbb L}^1$-norm  rather than the  ${\mathbb L}^{3/2}$-norm.

\section{Technical lemmas} \label{TL}

\setcounter{equation}{0}

Suppose that we have a sequence of random vectors $\{(A_j , B_j )\}_{1 \leq j \leq J} $ and a  filtration $\{{\mathcal H}_j  \}_{1 \leq j \leq J} $  such that 
$$\Big(\E_{{\mathcal H}_j} (A_j^2), \E_{{\mathcal H}_j} (|A_j|^p), \E_{{\mathcal H}_j} (B_j^2)\Big)_{j\in J}
$$ is a sequence of  independent random vectors (with values in ${\mathbb R}^3$). 
For any real $a$, let  
\[
H_j (a) = A_j+ a B_j  \, \text{ and  } \, \varphi_{j,a}^{\mathcal H} (x) = \E \big ( {\rm exp} ( {\rm i} x H_j(a) ) | {\mathcal H}_j  \big )  \, .
\]
With the notations above, the following modification of \cite[Lemma 4.5]{Ji} holds:

\begin{Lemma} \label{lma4.5} Let $p >2$. Let $J\ge 16$ be an integer. Assume the following:

(i) $\E_{{\mathcal H}_j} (A_j )=\E_{{\mathcal H}_j} (B_j )= 0 $, 
for any $1\le j\le J$,

(ii) there exists  $u^- >0$ such that ${\mathbb P} ( \E_{{\mathcal H}_j} (A^2_j  )  \leq u^- )< 1/2$, for any $1\le j\le J$,

(iii) $\sup_{j \geq 1} \E(|A_j|^p)  < \infty$,

(iv)  $\sup_{j \geq 1} \E(B_j^2) < \infty$.

\noindent Then there exist positive finite constants $c_1$, $c_2$ and $c_3$ depending only on $p$, $u^-$, $\sup_{j \geq 1} \E(|A_j|^p) $ and $\sup_{j \geq 1} \E(B_j^2) $ such that for any $a \in [0,c_1]$ and any 
$x^2 \leq c_2$, 
\[
 \E \Big  ( \prod_{j=1}^J |  \varphi_{j,a}^{\mathcal H} (x)|\Big )  \le {\rm e}^{-c_3 x^2 J} + {\rm e}^{- J/32}  \, . 
\]
\end{Lemma}

\noindent {\bf Proof of Lemma \ref{lma4.5}.} The beginning of the proof  proceeds as  the proof of \cite[Lemma 4.5]{Ji} but with substantial modifications.

\medskip

Let $1\le j\le J$ be fixed for the moment. Using a Taylor expansion we have
\[
 \E \big ( {\rm exp} ( {\rm i} x H_j(a) ) | {\mathcal H}_j  \big )  = 1 -   \E_{{\mathcal H}_j} (H^2_j (a))  x^2/2 + x^2/2 \int_0^1 (1-s) I(s,x) ds \, , 
\]
where, for any $h >0$ and any $s \in [0,1]$, 
\begin{align*}
 |  I(s,x) | & \leq 4 a^2 \E_{{\mathcal H}_j} (B^2_j )  + 2  \E_{{\mathcal H}_j} \big ( A^2_j   \big | ( \cos (sxH_j (a) ) - \cos (0) ) + {\rm i}  ( \sin (sxH_j (a) ) - \sin(0)  ) \big | \big )     \\
  & \leq 4 a^2 \E_{{\mathcal H}_j} (B^2_j )  + 8   \E_{{\mathcal H}_j} (A^2_j ) |xh| + 4  \E_{{\mathcal H}_j} (A^2_j  {\bf 1}_{|H_j (a) | \geq 2 h} ) \, .
\end{align*}
Using the fact that for any reals $u$ and $v$, $u^2 {\bf 1}_{|u+v| \geq 2 h}  \leq  u^2 {\bf 1}_{|u| \geq h} +  v^2$, we get 
\begin{align*}
|  I(s,x) | &  \leq  8 a^2 \E_{{\mathcal H}_j} (B^2_j ) + 8   \E_{{\mathcal H}_j} (A^2_j ) |xh| + 4  \E_{{\mathcal H}_j} (A^2_j  {\bf 1}_{|A_j  | \geq  h} ) \\
&  \leq  8 a^2 \E_{{\mathcal H}_j} (B^2_j ) + 8   \E_{{\mathcal H}_j} (A^2_j ) |xh| + 4 h^{2-p} \E_{{\mathcal H}_j} (|A_j|^{p}  )  \, .
\end{align*}
Now, for any $\alpha >0$,
\[
\big  \vert \E_{{\mathcal H}_j}(H^2_j (a)) - \big (  \E_{{\mathcal H}_j}(A_j^2
 ) + a^2 \E_{{\mathcal H}_j} (B^2_j )  \big )  \big  \vert  \leq \alpha^{-1}   \E_{{\mathcal H}_j}(A_j^2) + \alpha a^2 \E_{{\mathcal H}_j}( B^2_j) \, .
\]
So, overall, for any $h >0$ and any $\alpha >0$,
\begin{multline*}
 \Big |  \E \big ( {\rm exp} ( {\rm i} x H_j(a) ) | {\mathcal H}_j  \big )   -  1 +  \E_{{\mathcal H}_j} (A^2_j ))  x^2/2 \Big |  
 \leq   x^2 ( 3 a^2 +  \alpha a^2 ) \E_{{\mathcal H}_j} (B^2_j )  /2 \\  + 
 \E_{{\mathcal H}_j} (A^2_j ) ( x^2 \alpha^{-1} /2 + 2 h |x|^3 )   + x^2 h^{2-p} \E_{{\mathcal H}_j} (|A_j|^{p}  ) \, . 
\end{multline*}

Let us take $h = |x|^{-1/(p-1)} $.  Set $\delta (p) := (p-2)/(p-1)$.

\smallskip Let $\tilde u, u^+$ be positive numbers to be chosen later. 

\smallskip

Recall that by the conditional Jensen inequality, 
$\E_{{\mathcal H}_j} (A_j^2)\le \big(\E_{{\mathcal H}_j} (|A_j|^p)\big)^{2/p}$ ${\mathbb P}$-almost surely. For the sake of simplicity, we shall assume that this inequality takes place everywhere.

\smallskip

From the above computations, we infer that,  on the set $
\{\E_{{\mathcal H}_j} (B_j^2)\le \tilde u\} \cap \{\E_{{\mathcal H}_j} (|A_j|^p)\le u^+\}
$, one has, for any $\alpha >0$, 
\begin{multline*}
 \Big |  \E \big ( {\rm exp} ( {\rm i} x H_j(a) ) | {\mathcal H}_j  \big )   -  1 +  \E_{{\mathcal H}_j} (A^2_j ))  x^2/2 \Big |   \\
 \leq   x^2 ( 3 a^2 +   \alpha a^2 ) {\tilde u} /2 + 
x^2 (u^+)^{2/p} \alpha^{-1} /2 +  |x|^{2 + \delta(p)} ( 2  (u^+)^{2/p}   +  u^+ )  \, . 
\end{multline*}

Set
\[
u(x) : =   a^2  ( 3  +  \alpha) {\tilde u} /2+ 
 (u^+)^{2/p} \alpha^{-1} /2 +   |x|^{ \delta(p)} ( 2  (u^+)^{2/p}   +  u^+ ) \, .
\]

Let $u^-$ be a positive number ($u^-$ will be given by $(ii)$ but it is unimportant at this stage). We infer that, for every $x$ such that $x^2\le 2/u^-$ and $x^2\le 2/(u^+)^{2/p}$, on the set
$$
\Gamma_j:= \{\E_{{\mathcal H}_j} (B_j^2)\le \tilde u\}\cap 
\{\E_{{\mathcal H}_j} (A_j^2)> u^-\} \cap \{\E_{{\mathcal H}_j} (|A_j|^p)\le u^+\}\, 
$$
one has
$$
\Big |  \E \big ( {\rm exp} ( {\rm i} x H_j(a) ) | {\mathcal H}_j  \big ) \Big|\le 1-u^-x^2/2+x^2u(x)\, .
$$

Select now $\alpha = 8 (u^+)^{2/p}  / u^-$. Since $0 < u^-\! , u^+, {\tilde u} < \infty$, note that there exist positive 
constants $c_1, c_2 < \infty$ (depending only on  $( u^-, u^+, {\tilde u})$) such that 
\begin{gather*}
a \leq c_1  \Rightarrow a^2 ( 3  +  \alpha ) {\tilde u} /2 \leq u^-/16 \, ,\\ 
x^2\le c_2 \Rightarrow |x|^{ \delta(p)} ( 2  (u^+)^{2/p}   +  u^+ ) 
\le u^-/8\, .
\end{gather*}
Therefore,  there exist  constants $0 < c_1,c_2 < \infty$ (depending only on  $(\tilde u,  u^-\! , u^+)$) such that for any $a \leq  c_1$ and any $x^2 \leq c_2$, we have, on the set $\Gamma_j$,
\[
 \big |  \E \big ( {\rm exp} ( {\rm i} x H_j(a) ) | {\mathcal H}_j  \big )   \big | \le 1-u^-x^2/4 \leq e^{-  u^- x^2/4} \, .
\]

Set also $\Sigma_J:= \sum_{j=1}^J{\bf 1}_{\Gamma_j}$ and 
$\Lambda_J:= \{ \Sigma_J\ge J/8\}$.

\smallskip

From the previous computations and the trivial bound $\big |  \E \big ( {\rm exp} ( {\rm i} x H_j(a) ) | {\mathcal H}_j  \big )   \big | \le 1$, we see that, for any 
$0<\tilde u, u^-\! , u^+ <\infty$, there exist positive contants 
$c_1,c_2, c_3$ such that for every $x^2\le c_2$ and every $a\le c_1$, one has (recall that $J\ge 16$),
\[
\Big (  \prod_{j=1}^J |  \varphi_{j,a}^{\mathcal H} (x)| \Big )  {\bf 1}_{ \Lambda_J } \le 
{\rm e}^{-u^-x^2[J/8]/2}\le {\rm e}^{-u^-x^2J/32}\, .
\]

Using the above trivial bound again, the lemma will be proved if, with $u^-$ given by $(ii)$, one can chose 
$\tilde u, u^+>0$ such that $\BBP(\Lambda^c_J)\le {\rm e}^{-J/32}$.

\smallskip

By Markov's inequality and condition $(iv)$,
 \[\BBP(\E_{{\mathcal H}_j} (B_j^2)> \tilde u)
\le \frac{\sup_{j\in J}\E(B_j^2)}{\tilde u}
\underset{\tilde u\to +\infty}\longrightarrow 0\, .
\]
 Hence there exists 
$\tilde u>0$ such that, for any $1\le j\le J$,  $\BBP(\E_{{\mathcal H}_j} (B_j^2)> \tilde u)\le 1/8$.

Similarly, by condition $(iii)$, there exists $u^+>0$ such that, for any $1\le j\le J$,  $\BBP(\E_{{\mathcal H}_j} (|A_j|^p)> u^+)\le 1/8$.

On another hand,  by condition $(ii)$   and by definition of $\tilde u$ and $u^+$, we have 
\[
\E(\Sigma_J)\ge \sum_{j=1}^J(1- (1/2 + 1/8 + 1/8))=J/4\, .
\]

Hence,  
\begin{multline*}
\BBP(\Lambda^c_J)= \BBP (\Sigma_J <  J/8)  = \BBP (\Sigma_J  -  \E(\Sigma_J)<  J/8- \E(\Sigma_J))\\ \le 
\BBP (\Sigma_J  -  \E(\Sigma_J)<   - J/8) = \BBP ( -\Sigma_J  + \E(\Sigma_J)  >    J/8)  \, .
\end{multline*}
Therefore, using Hoeffding's inequality (see \cite[Theorem 2]{Hoeffding}),
\[
\BBP(\Lambda_J^c)\le {\rm e}^{\frac{-2(J/8)^2}{J}} =  {\rm e}^{-J/32}
\, ,
\]
which ends the proof of the lemma. \qed

\smallskip

For the next lemma, let us introduce the following notation: for any real $\beta$, let 
\beq \label{defkappalma}
 \kappa_\beta =   \frac{(\beta +1) (q-3/2)}{ q-1/2}  \, .
\eeq

\begin{Lemma} \label{lma4.7} Assume that $\mu$ has a moment of order $q  >2 $.   Let $X_{k,m}$ be defined by \eqref{defXkm}. Then, setting ${\bar X}_{k,m} = X_{k,m}  - \E_m (X_{k,m})$, for any real $\beta$ such that $-1 < \beta < q-3 + 1/q$, we have
\[
\Big \Vert \E_m \Big ( \sum_{k=m+1}^{2m} {\bar X}_{k,m}  \Big )^2  -  \E  \Big ( \sum_{k=m+1}^{2m} {\bar X}_{k,m}  \Big )^2  \Big \Vert_{1}\ll   1 + m^{3-q}  {\bf 1}_{q \leq 3}+  m^{ 1 -  \kappa_{\beta}} {\bf 1}_{ \beta < (q-3/2)^{-1}}   \, ,
\]
where $\kappa_\beta$ is defined in \eqref{defkappalma} and $\E_m (\cdot)$ means $\E(\cdot |  {\mathcal G}_m )$ with ${\mathcal G}_m = \sigma (W_0, \varepsilon_1, \ldots, \varepsilon_m)$.  In particular, if $q>3$, then
\[
\Big \Vert \E_m \Big ( \sum_{k=m+1}^{2m} {\bar X}_{k,m}  \Big )^2  -  \E  \Big ( \sum_{k=m+1}^{2m} {\bar X}_{k,m}  \Big )^2  \Big \Vert_{1}\ll  m^{1/5 }  \,  . 
\]
\end{Lemma}

\noindent{{\bf Proof of Lemma \ref{lma4.7}.} Note first that 
\begin{align} \label{step1lma47}
 & \Big \Vert \E_m \Big ( \sum_{k=m+1}^{2m} {\bar X}_{k,m}  \Big )^2    -  \E  \Big ( \sum_{k=m+1}^{2m} {\bar X}_{k,m}  \Big )^2  \Big \Vert_{1}  \\ &  \leq   \Big \Vert \E_m  \Big ( \sum_{k=m+1}^{2m}  X_{k,m}  \Big )^2 -   \E  \Big ( \sum_{k=m+1}^{2m}  X_{k,m}  \Big )^2 \Big  \Vert_{1}  + 2  \Big \Vert \E_m  \Big ( \sum_{k=m+1}^{2m}  X_{k,m}  \Big ) \Big  \Vert^2_{2}  \nonumber  \\
& := I_m + I\!\!I_m   \nonumber  \, .
\end{align}
Taking into account \eqref{Borne1condexpect},  \eqref{estimatedelta} and the fact that $q\geq 2$, we get 
\beq \label{P1lma47}
\Big \Vert \E_m  \Big (\sum_{k=m+1}^{2m}   X_{k,m}  \Big ) \Big  \Vert_{2} \ll  \sum_{k=m+1}^{2m}  \Vert \E_m   (   X_{k,m}   )  \Vert_{2}  \ll \sum_{k=1}^m \delta_{1, \infty} (k) \ll 1 \, .
\eeq
 It remains to handle $I_m$. With this aim, we first write the following decomposition: for any $\gamma \in (0,1]$
\begin{multline} \label{decIm}
 I_m \leq    \sum_{k=1}^{m} \Vert \E_m (  X^2_{k+m,m})-   \E ( X_{k+m,m} ^2 ) \Vert_{1}  \\ + 2  \sum_{\ell=1}^{m} \ell^{\gamma} \sup_{\ell \leq  j < i \leq \min ( 2 \ell  , m)} \Vert \E_m (  X_{i+m,m} X_{j+m,m})-   \E (  X_{i+m,m} X_{j+m,m})\Vert_{1}   \\
 + 2 \sum_{\ell=1}^{m} \sum_{k=[\ell^{\gamma} ]+1}^{m-\ell}\Vert \E_m (  X_{\ell +m,m} X_{\ell+k +m,m})-   \E (  X_{\ell+m,m} X_{\ell+k+m,m})\Vert_{1} \, .
\end{multline}
Note that for $1 \leq i,j \leq m$, 
\[
 \Vert \E_m (  X_{i+m,m} X_{j+m,m})-   \E (  X_{i+m,m} X_{j+m,m})\Vert_{1}  \leq \sup_{{\bar x}_1 , {\bar x}_2 \in X  \atop{{\bar y}_1 , {\bar y}_2 \in X}}  \E \big | X_{i, {\bar x}_1 }  X_{j, {\bar x}_2 } - X_{i, {\bar y}_1 }  X_{j, {\bar y}_2 } \big |  \, .
\]
With the same arguments as those developed in the proof of \cite[Prop. 4]{CDJ}, and  since $\mu $ has a moment of order $q>2$,   we then  infer  that 
\beq \label{pr4.5prime}
  \sum_{k \geq 1 } k^{q-3} \Vert \E_m (  X^2_{k+m,m})-   \E ( X_{k+m,m}^2 ) \Vert_{1}   \ll 1  \, , 
\eeq
and, for every $\beta <  q-3+1/q$, 
\beq \label{item2prop4}
 \sum_{\ell \geq  1} \ell^{ \beta} \sup_{\ell \leq  j < i \leq \min (2 \ell,m) } \Vert \E_m (  X_{i+m,m} X_{j+m,m})-   \E (  X_{i+m,m} X_{j+m,m})\Vert_{1}  \ll 1  \, . 
\eeq
On another hand, with the same arguments as those used to prove \cite[Relation (34)]{CDJ}, we first write 
\begin{multline*}
\sum_{\ell=1}^{m} \sum_{k=[\ell^{\gamma} ]+1}^{m-\ell}\Vert \E_m (  X_{\ell +m,m} X_{\ell+k +m,m})-   \E (  X_{\ell+m,m} X_{\ell+k+m,m})\Vert_{1} \\
\ll \Big ( \sum_{\ell=m+1}^{2m}  \Vert \E_m(X_{\ell,m} )\Vert_{2} \Big )^2 + \sum_{\ell=1}^{m} \sum_{k=[\ell^{\gamma} ]+1}^{m-\ell} \sum_{u=1}^{\ell}  \Vert P_{m+1} ( X_{u + m,m} ) \Vert_2 \Vert P_{m+1} ( X_{u +k+ m,m} ) \Vert_2 \, ,\\
\ll \Big ( \sum_{\ell=m+1}^{2m}  \Vert \E_m(X_{\ell,m} )\Vert_{2} \Big )^2  + \Big(\sum_{v = 1}^m a(0,v) \Big)\Big(
\sup_{u\ge 1} \sum_{\ell=1}^{m} \sum_{k=[\ell^{\gamma} ]+1}^{m-\ell}  a(k,u) \, \Big) \, , 
\end{multline*}
where we have used the notations $P_{m+1}( \cdot) =  \E_{m+1}  ( \cdot)  -  \E_{m}( \cdot) $ and  $a(k,u) = \Vert P_{m+1} ( X_{u +k+ m,m} ) \Vert_2$. Note first that  
\begin{align*}
 \sum_{\ell=1}^{m} \sum_{k=[\ell^{\gamma} ]+1}^{m-\ell}   a(k,u) &  \ll  \sum_{k=2}^{m-1}  ( k^{1/\gamma} \wedge m )  a(k,u)  \\&  \ll
    \sum_{k=2}^{[m^{\gamma} ] }  k^{-1} a(k,u)  \sum_{ \ell =1}^k \ell^{1 / \gamma}  + m   \sum_{k =  [m^{\gamma} ]  +1}^m  k^{-1}a(k,u)  \sum_{ \ell =1}^k 1  \, . \end{align*}
Changing the order of summation and using Cauchy-Schwarz's inequality, it follows that 
\begin{multline*}
 \sum_{\ell=1}^{m} \sum_{k=[\ell^{\gamma} ]+1}^{m-\ell}   a(k,u) 
 \ll  \sum_{\ell =1}^{[m^{\gamma} ] }   \ell^{1 / \gamma - 1/2}  \Big (  \sum_{ k \geq \ell } a^2(k,u)  \Big )^{1/2}  \\ + m \sum_{\ell =   [m^{\gamma}  ]+1 }^{m }    \ell^{ - 1/2}  \Big (  \sum_{ k \geq  \ell }  a^2(k,u)    \Big )^{1/2}  +    m^{ 1 + \gamma/2}  \Big (  \sum_{ k \geq  [m^{\gamma}  ]+1 }  a^2(k,u)   \Big )^{1/2}   \, .
\end{multline*}
But, for any $u \geq 1$, by stationarity, 
\[
 \Big (  \sum_{ k \geq  \ell }  a^2(k,u)    \Big )^{1/2}   \leq   \Vert \E_{m+1} ( X_{u +\ell+ m,m} ) \Vert_2 \leq  \Vert \E_{m+1} ( X_{u +\ell+ m,m} ) \Vert_{\infty} \leq \delta_{1, \infty} ( \ell)  \, .
\]
Notice also that 
\[
\sum_{v= 1}^m a(0,v)\le \sum_{v= 1}^m
\Vert \E_{m+1}(X_{v+m,m})\Vert_2 \le \sum_{v= 1}^m
\delta_{1,\infty}(v) \, .
\]

\smallskip

Hence, from the above considerations and taking into account \eqref{P1lma47}, \eqref{estimatedelta}, \eqref{estimatedeltabis} and the fact that $\mu$ has a moment of order $q \geq 2$, we infer that 
\begin{multline}  \label{pr4.7}
\sum_{\ell=1}^{m} \sum_{k=[\ell^{\gamma} ]+1}^{m-\ell}\Vert \E_m (  X_{\ell +m,m} X_{\ell+k +m,m})-   \E (  X_{\ell+m,m} X_{\ell+k+m,m})\Vert_{1} \\
\ll 1  + m^{1 - \gamma  ( q -  3/2 ) }  {\bf 1}_{1 /\gamma > q-3/2 }  \, .
\end{multline}
Starting from \eqref{decIm} and considering the estimates \eqref{pr4.5prime}, \eqref{item2prop4} and  \eqref{pr4.7}, 
we get, for any $\gamma \in (0,1]$ and any $\beta$ such that $-1 < \beta < q-3 + 1/q$, 
\beq \label{Imlma47new}
I_m  \ll    1 + m^{3-q}  {\bf 1}_{q \leq 3}+  m^{\gamma- \beta} {\bf 1}_{\gamma > \beta} +m^{1 - \gamma  ( q -  3/2 ) }  {\bf 1}_{1 /\gamma > q-3/2 }   \, .
\eeq
Let us select now $ \gamma$ such that $\gamma - \beta = 1 - \gamma  ( q -  3/2 )$.  This gives $ \gamma = ( \beta +1) / (q-1/2)$. Since $\beta >-1$, $\beta < q-3 + 1/q$ and $q >2$ we have  $\gamma  \in (0, 1]$. Moreover $1 /\gamma > q-3/2 $ and $\gamma > \beta$ provided $ \beta < (q-3/2)^{-1}$.  Starting from \eqref{step1lma47} and taking into account \eqref{P1lma47}, \eqref{Imlma47new}  and the above selection of $\gamma$, which entails that $\kappa_\beta = \gamma ( q-3/2)$,  the lemma follows. \qed

\begin{Lemma} \label{lmaR1normep} Let $p  \geq 2$. Assume that $\mu$ has a moment of order $q$  in $ ]p, p+1]$.  Then $ \Vert R_1 \Vert^p_{p} \ll m^{p+1-q}$, where $R_1$ is defined by \eqref{defRj}.   
 \end{Lemma}

\noindent{{\bf Proof of Lemma \ref{lmaR1normep}.}  Let ${\tilde X}_{k,m} = X_{k,m} - \E_{{\mathbb F}_m} (X_{k,m})$ and 
$\E_\ell (\cdot):=\E(\cdot |  {\mathcal G}_\ell )$ with ${\mathcal G}_\ell = \sigma (W_0, \varepsilon_1, \ldots, \varepsilon_\ell)$.   We write 
\[
{\tilde X}_{k,m}  = ({\tilde X}_{k,m}  - \E_{k-1} ( {\tilde X}_{k,m}  )) +  \E_{k-1} ( {\tilde X}_{k,m}  ) := d_{k,m} + r_{k,m} \, ,
\]
and then 
\begin{equation} \label{dec1R1mart}
\Vert R_1 \Vert_{p} \leq \Big  \Vert  \sum_{k=m+1}^{2m}   d_{k,m}  \Big \Vert_p +  \Big  \Vert  \sum_{k=m+1}^{2m}   r_{k,m}  \Big \Vert_p \, .
\end{equation}
Note that $(d_{k,m})_{k \geq 1} $ is a sequence of ${\mathbb L}^q$-martingale differences with respect to the filtration $({\mathcal G}_{k})_{k \geq 1} $. Moreover, for any $r \geq 1$, 
$\Vert d_{k,m} \Vert_r \leq 2  \Vert {\tilde X}_{k,m} \Vert_r$ and, for any integer $k \in [m+1, 2m]$, 
  \begin{multline*} 
 \E \big | {\tilde X}_{k,m}  \big |^r \\
   = \E \Big | f_m (\varepsilon_{k-m+1}, \ldots, \varepsilon_{m}, \varepsilon_{m+1}, \ldots,  \varepsilon_{k})  - \int   f_m (v_{k-m+1}, \ldots, v_{m}, \varepsilon_{m+1}, \ldots,  \varepsilon_{k})  \prod_{i=k-m+1}^{m} d\mu (v_i) \Big |^r \\
   \leq \int  \E \Big | f_m (\varepsilon_{k-m+1}, \ldots, \varepsilon_{m}, \varepsilon_{m+1}, \ldots,  \varepsilon_{k})  -   f_m (v_{k-m+1}, \ldots, v_{m}, \varepsilon_{m+1}, \ldots,  \varepsilon_{k}) \Big |^r   \prod_{i=k-m+1}^{m} d\mu (v_i)   
 \, .
\end{multline*}
Hence, for any integer $k \in [m+1, 2m]$ and any $r \geq 1$, 
 \begin{align} \label{P1lmaR1norme2}
 \E \big | {\tilde X}_{k,m}  \big |^r 
&  \leq   \int  \!\! \int \E \Big | f_m (u_{k-m+1}, \ldots, u_{m}, \varepsilon_{m+1}, \ldots,  \varepsilon_{k})  \nonumber  \\& \quad \quad  -   f_m (v_{k-m+1}, \ldots, v_{m}, \varepsilon_{m+1}, \ldots,  \varepsilon_{k}) \Big |^r   \prod_{i=k-m+1}^{m} d\mu (v_i)    \prod_{i=k-m+1}^{m} d\mu (u_i)    \nonumber \\
&  \leq \sup_{{\bar x} , {\bar y} \in X} \E  | X_{k-m,{\bar x}} - X_{k-m,{\bar y}}|^r= \delta^r_{r,\infty}\ (k-m) \,  .
\end{align}
On another hand $(r_{k,m})_{k \geq 1} $ is a sequence of centered random variables  such that 
\[
\Vert r_{k,m} \Vert_{\infty} \leq  2 \Vert \E ( |X_k| | {\mathcal G}_{k-1} ) \Vert_{\infty} \leq 2   \int_G \log (N(g)) \mu (dg ):= K < \infty \, .
\]

To handle the first term in the right-hand side of \eqref{dec1R1mart}, we use the Rosenthal-Burkholder's inequality for martingales (see \cite{Bu73}). Hence, there exists a positive constant $c_p$ only depending on $p$ such that 
\begin{equation*}
 \Big \Vert \sum_{k=m+1}^{2m}   d_{k,m}  \Big \Vert_p^p  \le c_p\Big \{ \Big \Vert \sum_{k=m+1}^{2m}  \E(d_{k,m}^2| {\mathcal G}_{k-1}) \Big \Vert_{p/2}^{p/2} +\sum_{k=m+1}^{2m}  \Vert d_{k,m } \Vert_p^p\Big \} \, .
\end{equation*}
Taking into account \eqref{P1lmaR1norme2},  \eqref{estimatedelta} and the fact that $\mu$ has a moment of order $q=p+1$, it follows that 
\begin{equation*}
\sum_{k=m+1}^{2m}  \Vert d_{k,m } \Vert_p^p \leq 2^p  \sum_{k=m+1}^{2m}    \delta^p_{p,\infty}\ (k-m) \ll m^{p+1-q}\, .
\end{equation*}
On another hand, by the properties of the conditional expectation, note that 
\[
 \Vert  \E(d_{k,m}^2| {\mathcal G}_{k-1}) \Vert_\infty \leq  \Vert \E ( X^2_k | {\mathcal G}_{k-1} ) \Vert_{\infty} \leq  \int_G  ( \log (N(g)) )^2 \mu (dg ):= L < \infty \, .
\]
Hence, by using  \eqref{P1lmaR1norme2}, 
\[
\Vert \E(d_{k,m}^2| {\mathcal G}_{k-1})  \Vert_{p/2}^{p/2} \le  L^{(p-2)/2}  \Vert  d_{k,m} \Vert_2^2 \leq  4L^{(p-2)/2}  \Vert  {\tilde X}_{k,m} \Vert_2^2 \leq   4L^{(p-2)/2}   \delta^2_{2,\infty}\ (k-m)  \, .
\]
It follows that 
\[ 
\Big \Vert \sum_{k=m+1}^{2m}  \E(d_{k,m}^2| {\mathcal G}_{k-1}) \Big \Vert_{p/2}^{p/2} 
\le  4  L^{(p-2)/2}  
 \Big ( \sum_{k= 1}^m  \delta^{4/p}_{2,\infty} (k)  \Big )^{p/2} \, .
\]
By taking into account  \eqref{estimatedelta} (when $p=2$) and \eqref{estimatedeltabis}  (when $p >2$), and since $q \in ]p, p+1]$, we get 
\[ 
\Big \Vert \sum_{k=m+1}^{2m}  \E(d_{k,m}^2| {\mathcal G}_{k-1}) \Big \Vert_{p/2}^{p/2} 
\le  4  L^{(p-2)/2}   m^{p+1-q} \, .
\]
So, overall, 
\begin{equation}\label{burkholder-2}
 \Big \Vert \sum_{k=m+1}^{2m}   d_{k,m}  \Big \Vert_p^p \ll m^{p+1-q} \, .
\end{equation}

We handle now the second term  in the right-hand side of \eqref{dec1R1mart}.  By using the Burkholder-type inequality stated in \cite[Proposition 4]{DD03}, we get
\[
 \Big \Vert \sum_{k=m+1}^{2m}   r_{k,m}  \Big \Vert_p^2 \leq 2 p  \sum_{i= m+1}^{2m} \sum_{k=i}^{2m} \Vert   r_{i,m} \E (  r_{k,m} | {\mathcal G}_{i-1}  ) \Vert_{p/2} \, .
\]
For any $ k \geq i$, by the computations leading to the upper bound \cite[(63)]{CDM}, we have 
\beq \label{condexpectationrkm}
 \Vert  \E (  r_{k,m} | {\mathcal G}_{i-1}  ) \Vert_{\infty} \leq  \delta_{1, \infty} (k-i +1) \, , 
\eeq
implying that 
\[
 \Vert   r_{i,m} \E (  r_{k,m} | {\mathcal G}_{i-1}  ) \Vert_{p/2} \ \leq  \Vert   r_{i,m}   \Vert_{p/2}  \delta_{1, \infty} (k-i +1) \, .
\]
Since $\mu$ has a moment of order at least $2$, by \eqref{estimatedelta}, $\sum_{\ell \geq 1}   \delta_{1, \infty} (\ell) < \infty$. Hence 
\[
 \Big \Vert \sum_{k=m+1}^{2m}   r_{k,m}  \Big \Vert_p^2 \ll \sum_{i= m+1}^{2m}    \Vert   r_{i,m}  \Vert_{p/2} \, .
\]
But, for any $r \geq 1$,   $\Vert   r_{i,m}  \Vert^r_{r}  \leq K^{r-1}  \Vert   r_{i,m}  \Vert_1  \leq 2 K^{r-1}  \Vert  {\tilde X}_{i,m} \Vert_1 $. Hence, by using \eqref{P1lmaR1norme2}, it follows that, 
for any $r \geq 1$, $\Vert   r_{i,m}  \Vert^r_{r}   \leq 2 K^{r-1}  \delta_{1, \infty} (i-m ) $. Therefore, by  \eqref{estimatedeltabis} and the fact that $q-1 >p/2$ (since $q >p$ and  $p \geq 2$), we derive that 
\begin{equation}\label{burkholder-3}
 \Big \Vert \sum_{k=m+1}^{2m}   r_{k,m}  \Big \Vert_p^p
\ll   \Big ( 
\sum_{i= 1}^{m}   \delta^{2/p}_{1, \infty} (i ) \Big )^{p/2}  \ll 1  \, .
\end{equation}
Starting from \eqref{dec1R1mart} and considering the upper bounds \eqref{burkholder-2} and \eqref{burkholder-3}, the lemma follows. \qed
  
\begin{Lemma} \label{lmamomentp} Assume that $\mu$ has a finite moment of order $q \geq 2$. Then $
\big \Vert \sum_{k=m+1}^{2m } X_{k} \big  \Vert_q \ll \sqrt{m}
$ and 
$
\big \Vert \sum_{k=m+1}^{2m } X_{k,m} \big  \Vert_q \ll \sqrt{m}
$. 
\end{Lemma}

\noindent{{\bf Proof of Lemma \ref{lmamomentp}}.}  The two upper bounds are proved similarly. Let us prove the second one. As to get \eqref{momentAJ}, we use \cite[Cor. 3.7]{MPU19}, to derive that 
\[
\Big \Vert \sum_{k=m+1}^{2m } X_{k,m} \Big  \Vert_q  \ll  \sqrt{m} \Big [  \Vert X_{1+m,m} \Vert_{q} + \sum_{k=m+1}^{2m} k^{-1/2} \Vert \E_m (X_{k,m}) \Vert_{q}  \Big ]  \, ,
\]
where $\E_m (\cdot)$ means $\E(\cdot |  {\mathcal G}_m )$ with ${\mathcal G}_m = \sigma (W_0, \varepsilon_1, \ldots, \varepsilon_m)$. 
But $\Vert X_{1+m,m} \Vert_{q}  \leq \Vert X_{1} \Vert_{q}  < \infty$ and  $\Vert \E_m (X_{k+m,m}) \Vert_{q} \leq \Vert \E_m (X_{k+m,m}) \Vert_{\infty}\leq \delta_{1, \infty} (k)$. Hence, the lemma follows by considering   \eqref{estimatedelta}. \qed 

\medskip

For the next lemma, we recall the notations  \eqref{defbbFm}  and  \eqref{notaSm1} for ${\mathbb F}_m$ and $Y_j^{(1)}$.

 \begin{Lemma} \label{lma4.9} Assume that $\mu$ has a finite moment of order $q \in ]2,3]$. Then for $f(x) \in \{ \cos x, \sin x \}$, we have 
\[
\Big \Vert \E_{{\mathbb F}_m}  \Big [ f \Big ( t \frac{Y_2^{(1)}}{\sqrt{2m}} \Big ) \Big ]  -  \E \big [ f  ( t  s N /{\sqrt 2} ) \big ] \Big   \Vert_{1}  \ll  \frac{ t^2 }{m^{q/2 -1}}  + \frac{|t|}{m^{q-3/2}}  \, .
\]
In addition
\[
\Big \Vert \E_{{\mathbb F}_m}  \Big [ f \Big ( t \frac{Y_1^{(1)}}{\sqrt{2m}} \Big ) \Big ]  -  \E \big [ f  ( t  s N /{\sqrt 2} ) \big ] \Big   \Vert_{1}  \ll  \frac{ t^2 }{m^{q/2 -1}} \, .
\]
\end{Lemma}

\noindent{{\bf Proof of Lemma \ref{lma4.9}.}  Since the derivative of $x \mapsto f(tx)$ is $t^2$-Lipschitz, making use of a Taylor expansion  as done in the proof of Item (2) of \cite[ Lemma 5.2]{DMR09}, we have 
\begin{multline} \label{P1lma49}
\Big \Vert \E_{{\mathbb F}_m}  \Big [ f \Big ( t \frac{Y_2^{(1)}}{\sqrt{2m}} \Big ) \Big ]  -  \E \big [ f  ( t  s N /{\sqrt 2} ) \big ] \Big   \Vert_{1} \\
\leq  \Big \Vert \E_{{\mathbb F}_m}  \Big [ f \Big ( t \frac{U_2}{\sqrt{2m}} \Big ) \Big ]  -  \E \big [ f  ( t  s N /{\sqrt 2} ) \big ] \Big   \Vert_{1} + \frac{ t^2 }{ 2m} \big (  \Vert R_2 \Vert_{2}   \Vert U_2 \Vert_{2}  +  \Vert R_2 \Vert^2_{2}   \big )   \, .
\end{multline}
Now recall that $U_2= \sum_{k=2m+1}^{3m}  {\tilde X}_{k,m}$ where ${\tilde X}_{k,m}=  X_{k,m} - \E_{{\mathbb F}_m} (X_{k,m} ) $ with 
 $X_{k,m} = \E ( X_{k} | {\mathcal E}_{k-m+1}^{k} ) :=f_m (  \varepsilon_{k-m+1}, \ldots, \varepsilon_{k})$. 
Let $(\varepsilon^*_k)_k$ be an independent copy of $(\varepsilon_k)_k$ and independent of $W_0$. Define
\beq \label{defU2*}
{X}^*_{k,m}= f_{m} (\varepsilon^*_{k-m+1}, \ldots,\varepsilon_{2m}^*, \varepsilon_{2m+1}, \ldots,  \varepsilon_{k})  \mbox{ and } U_2^*= \sum_{k=2m+1}^{3m} { X}^*_{k,m} \, .
\eeq
 Clearly $U_2^*$ is independent of $ {\mathbb F}_m $. Using again the fact that the derivative of $x \mapsto f(tx)$ is $t^2$-Lipschitz and a Taylor expansion  as  in the proof of \cite[ Lemma 5.2]{DMR09}, we get 
\begin{multline} \label{P2lma49}
 \Big \Vert \E_{{\mathbb F}_m}  \Big [ f \Big ( t \frac{U_2}{\sqrt{2m}} \Big ) \Big ]  - \E  \big [ f  ( t  s N /{\sqrt 2} ) \big ] \Big   \Vert_{1}  \\ \ll  
  \Big |\E  \Big [ f \Big ( t \frac{U^*_2}{\sqrt{2m}} \Big ) \Big ]  -  \E \big [ f  ( t s N  /\sqrt{2}) \big ] \Big   |  +  \frac{ t^2 }{2 m }   \big (  \Vert U_2 - U_2^* \Vert_{2}  \Vert U_2^* \Vert_{2}  + \Vert U_2 - U_2^* \Vert^2_{2}  \big )  \, .
\end{multline}
Setting 
${\mathcal G}_{k,m} = \sigma ( \varepsilon^*_{m+2}, \ldots,\varepsilon_{2m}^*, \varepsilon_{m+2}, \ldots,\varepsilon_{2m}, \varepsilon_{2m+1}, \ldots,  \varepsilon_{k}) $, we have
\begin{multline*} 
 \Vert U_2 - U_2^* \Vert^2_{2}  \leq  2 \Big (   \sum_{k=2m+1}^{3m} \Vert  \E_{{\mathbb F}_m} (X_{k,m} )  \Vert_{2} \Big )^2  +2   \sum_{k=2m+1}^{3m} \Vert   { X}_{k,m} -  { X}^*_{k,m}\Vert_{2}^2  \\
  + 4 \sum_{k=2m+1}^{3m} \sum_{\ell = k+1}^{3m}  \Vert ( {X}_{k,m} -  { X}^*_{k,m}) \E ( { X}_{\ell,m} -  { X}^*_{\ell,m}  | {\mathcal G}_{k,m})  \Vert_1 \, .
\end{multline*}
But, for any integer $k$ in $[2m+1,3m]$ and any $r \geq 1$,
\beq  \label{new1EFmr}
\Vert  \E_{{\mathbb F}_m} (X_{k,m} )  \Vert_{2} \leq \Vert  \E(X_{k,m}  | {\mathcal G}_{2m})  \Vert_{r}   \leq   \Vert  \E(X_{k,m}  | {\mathcal G}_{2m})  \Vert_{\infty} \leq \delta_{1, \infty} (k-2m) \, ,
\eeq
where ${\mathcal G}_{2m} = \sigma (W_0, \varepsilon_1, \ldots, \varepsilon_{2m})$. 
On another hand, proceeding as in the proof of \cite[Lemma 24]{CDM}, we get that, for any $k \geq 2m$ and  any $r \geq 1$, 
\beq  \label{new2EFmr}
\Vert    { X}_{k,m} -  { X}^*_{k,m} \Vert^r_{r}  \leq \delta_{r, \infty}^r (k-2m) \, . 
\eeq
Let us now handle the quantity $\Vert \E (  { X}_{\ell,m} -  { X}^*_{\ell,m}  | {\mathcal G}_{k,m})  \Vert_{\infty} $ for  $\ell >k$.  For this aim, let  $(\varepsilon'_k)_k$ be an independent copy of $(\varepsilon_k)_k$, independent also of $ ( (\varepsilon^*_k)_k,W_0)$.  With the notation ${\mathcal H}_{k,m }= \sigma ( (\varepsilon_{i})_{i\leq k},  W_0, \varepsilon^*_{m+1}, \ldots, \varepsilon^*_{2m} )$, one  has, for any integers $k, \ell$  in $ [2m+1, 3m]$ such that $\ell >k$, 
\begin{multline*} 
 \E (  { X}_{\ell,m} -  { X}^*_{\ell,m}  | {\mathcal G}_{k,m})  = \E  \big (f_{m} (\varepsilon_{\ell-m+1}, \ldots,\varepsilon_{2m}, \varepsilon_{2m+1}, \ldots,  \varepsilon_{k}, \varepsilon'_{k+1}  \ldots, \varepsilon'_{\ell})   | {\mathcal H}_{k,m }  \big )   \\
 -\E  \big (f_{m} (\varepsilon^*_{\ell-m+1}, \ldots,\varepsilon_{2m}^*, \varepsilon_{2m+1}, \ldots,  \varepsilon_{k}, \varepsilon'_{k+1}, \ldots, \varepsilon'_{\ell})   | {\mathcal H}_{k,m }  \big ) \, .
\end{multline*}
Therefore, by simple arguments and using stationarity, we infer  that, for $ k, \ell$  in $ [2m+1, 3m]$ such that $\ell >k$, 
\beq  \label{new3EFmr}
\Vert \E (  { X}_{\ell,m} -  { X}^*_{\ell,m}  | {\mathcal G}_{k,m})  \Vert_{\infty} \leq \sup_{{\bar x}, {\bar y} \in X} |\E (X_{ \ell - k, {\bar x}} )  - \E (X_{\ell -k , {\bar y}} )| \leq  \delta_{1, \infty} (\ell-k)  \, . 
\eeq
So, overall,
\[
 \Vert U_2 - U_2^* \Vert^2_{2}  \ll   \sum_{k=1}^m \delta_{2, \infty}^2 (k) +  \Big ( \sum_{k=1}^m \delta_{1, \infty} (k)  \Big )^2 \, .
\]
Taking into account \eqref{estimatedelta} and the fact that $\mu$ has a moment of order $q \in ]2,3]$,  it follows that 
\begin{equation}  \label{P4lma49}
 \Vert U_2 - U_2^* \Vert^2_{2}  \ll    m^{3-q} \, .
\end{equation}
On another hand, by stationarity, $ \Vert R_2 \Vert_{2} =  \Vert R_1 \Vert_{2}$, and by Lemma \ref{lmaR1normep}, since $\mu$ has a moment of order $q\in ]2,3]$, we have  $ \Vert R_1 \Vert_{2}   \ll m^{(3-q)/2}$.} Moreover, by using \eqref{P4lma49}, Lemma \ref{lmamomentp} 
and the fact that ${ X}^*_{k,m}$ is distributed as $X_{k,m}$, we get that  $\Vert U_2 \Vert_{2} + \Vert U^*_2 \Vert_{2}  \ll \sqrt{m}$. 
So, the inequalities \eqref{P1lma49}, \eqref{P2lma49} and \eqref{P4lma49} together with the above considerations, lead to 
\begin{equation} \label{P3lma49}
\Big \Vert \E_{{\mathbb F}_m}  \Big [ f \Big ( t \frac{Y_2^{(1)}}{\sqrt{2m}} \Big ) \Big ]  - \E  \big [ f  ( t  s N /{\sqrt 2} ) \big ] \Big   \Vert_{1}   \ll  
  \Big |\E  \Big [ f \Big ( t \frac{U^*_2}{\sqrt{2m}} \Big ) \Big ]  -  \E \big [ f  ( t  s N  /\sqrt{2}) \big ] \Big   |  +  \frac{ t^2 }{m^{q/2 -1}}    \, .
\end{equation}
Next,  taking into account that $x \mapsto f(tx)$ is $t$-Lipschitz and the fact that  $U_2^*=^{\mathcal D} \sum_{k=1}^{m}  X_{k+m,m}$ and $S_m=^{\mathcal D}S_{2m}-S_m$, we get 
\[
  \Big | \E  \Big [ f \Big ( t \frac{U^*_2}{\sqrt{2m}} \Big ) \Big ]  -  \E  \Big [ f \Big ( t \frac{S_{m}}{\sqrt{2m}} \Big ) \Big ]  \Big |  \leq \frac{|t|}{\sqrt{2m}}  \Big \Vert \sum_{k=1}^m ( X_{k+m,m} - X_{k+m}   ) \Big \Vert_{1} \, .
\]
But,  by stationarity, \cite[Lemma 24]{CDM} and \eqref{estimatedelta},  we have
\[
 \Big \Vert \sum_{k=1}^m ( X_{k+m,m} - X_{k+m}  ) \Big \Vert_{1} \leq m  \delta_{1, \infty} (m)  \ll 1/m^{q-2} \, ,
\]
implying that 
\beq  \label{P5lma49}
  \Big | \E  \Big [ f \Big ( t \frac{U^*_2}{\sqrt{2m}} \Big ) \Big ]  -  \E  \Big [ f \Big ( t \frac{S_{m}}{\sqrt{2m}} \Big ) \Big ]  \Big |  \ll \frac{|t|}{m^{q-3/2}} \, .
\eeq
Hence starting from \eqref{P3lma49} and taking into account  \eqref{P5lma49}, we derive that 
\begin{multline} \label{P6lma49}
\Big \Vert \E_{{\mathbb F}_m}  \Big [ f \Big ( t \frac{Y_2^{(1)}}{\sqrt{2m}} \Big ) \Big ]  -  \E  \big [ f  ( t sN / \sqrt{2} ) \big ] \Big   \Vert_{1} 
\\ \ll   \Big |\E  \Big [ f \Big ( t \frac{S_m}{\sqrt{2m}} \Big ) \Big ]  -  \E \big [ f  ( t s  N /\sqrt{2}) \big ] \Big   |  +\frac{ t^2 }{m^{q/2 -1}}   + \frac{ |t| }{m^{q-3/2}}   \, .
\end{multline}
Next note that $ x \mapsto f(tx)$ is such that its first derivative is $t^2$-Lipshitz. Hence, by the definition of the Zolotarev distance of order $2$ (see for instance the introduction of \cite{DMR09}  for the definition of those distances), 
\[
 \Big |\E  \Big [ f \Big ( t \frac{S_m}{\sqrt{2m}} \Big ) \Big ]  -  \E \big [ f  ( t  s N  /{\sqrt 2}) \big ] \Big   | \leq t^2 \zeta_2 \big (  P_{S_m/{\sqrt {2  m}}} , G_{s^2/2} \big) \, .
\]
We apply \cite[Theorems 3.1 and 3.2]{DMR09} and, since $\mu$ has a finite moment of order $q \in ]2,3]$, we derive
\[
 \zeta_2 \big (  P_{S_m/{\sqrt{2 m}}} , G_{s^2/2} \big) \ll m^{- ( q/2-1)} \, .
\]
Note that the fact that the conditions (3.1), (3.2), (3.4) and (3.5) required in  \cite[Theorems 3.1 and  3.2]{DMR09}  hold when $\mu$ has a finite moment of order $q \in ]2,3]$ has been established in the proof of  \cite[Theorem 2]{CDJ}. 
Hence
\begin{equation} \label{P7lma49}
 \Big |\E  \Big [ f \Big ( t \frac{S_m}{\sqrt{2m}} \Big ) \Big ]  -  \E \big [ f  ( t  s N /\sqrt{2}) \big ] \Big   | \ll \frac{ t^2 }{m^{q/2 -1}}   \, .
\end{equation}
Starting from \eqref{P6lma49} and considering  \eqref{P7lma49}, the first part of Lemma \ref{lma4.9} follows.  Now to prove the second part, we note that 
\begin{multline*} 
\Big \Vert \E_{{\mathbb F}_m}  \Big [ f \Big ( t \frac{Y_1^{(1)}}{\sqrt{2m}} \Big ) \Big ]  -  \E \big [ f  ( t  s N /{\sqrt 2} ) \big ] \Big   \Vert_{1} \\
\leq  \Big \Vert \E \Big [ f \Big ( t \frac{S_m}{\sqrt{2m}} \Big ) \Big ]  -  \E \big [ f  ( t  s N /{\sqrt 2} ) \big ] \Big   \Vert_{1} + \frac{ t^2 }{ 2m} \big (  \Vert R_1 \Vert_{2}   \Vert S_m \Vert_{2}  +  \Vert R_1 \Vert^2_{2}   \big )   \, ,
\end{multline*}
where we used the fact that $S_m$ is independent of ${\mathbb F}_m$. Hence the second part of Lemma \ref{lma4.9} follows by using  \eqref{P7lma49},  Lemma \ref{lmaR1normep} and the fact that, by Lemma \ref{lmamomentp},  $\Vert S_m \Vert_2 \ll \sqrt{m}$. \qed

 \begin{Lemma} \label{lmaU1U1starnormep}  Let $p  \geq 2$. Assume that $\mu$ has a moment of order $q$  in $ ]p, p+1]$. Then $ \Vert U_2 - U_2^* \Vert^p_{p} \ll m^{p+1-q}$, where $U_2$ is defined by \eqref{defUj} and $U^*_2$ is defined by \eqref{defU2*}.
 \end{Lemma}

\noindent{{\bf Proof of Lemma \ref{lmaU1U1starnormep}.}   When $p=2$, the lemma has been proved in \eqref{P4lma49}. Let us complete the proof for any $p  \geq 2$.  We shall follow the same strategy as in the proof of Lemma   \ref{lmaR1normep}. Let 
$Z_{k,m}:= X_{k,m} - { X}^*_{k,m}$ where ${X}^*_{k,m}$ is defined  by \eqref{defU2*}. Setting ${\mathcal F}^Z_j = \sigma (  \varepsilon_{m+2}, \ldots, \varepsilon_{j},   \varepsilon^*_{m+2}, \ldots, \varepsilon^*_{2m})$, 
\[
d_{k,m}^Z :=  Z_{k,m} - \E ( Z_{k,m} | {\mathcal F}^Z_{k-1} ) \, \text{ and } \,  r_{k,m}^Z = \E ( Z_{k,m} | {\mathcal F}^Z_{k-1} ) \, , 
\]
we have 
\begin{equation} \label{dec1UUstarmart}
\Vert U_2 - U_2^*  \Vert_{p} \leq   \sum_{k=2m+1}^{3m} \Vert \E ( X_{k,m} | {\mathbb F}_m  )\Vert_{p}   + \Big  \Vert  \sum_{k=2m+1}^{3m}   d^Z_{k,m}  \Big \Vert_p +  \Big  \Vert  \sum_{k=2m+1}^{3m}   r^Z_{k,m}  \Big \Vert_p \, .
\end{equation}
Recall the notation $  {\mathcal G}_{\ell} = \sigma (W_0, \varepsilon_1, \ldots, \varepsilon_{\ell})$. Note that 
\begin{equation} \label{dec1UUstarmart-1}
 \Vert \E ( ( d^Z_{k,m}  )^2  | {\mathcal F}^Z_{k-1} ) \Vert_{\infty} \leq  4  \Vert \E ( X^2_{k,m}   | {\mathcal G}_{k-1} ) \Vert_{\infty}  \leq 4   \int_G  ( \log (N(g)) )^2 \mu (dg )  < \infty
\end{equation}
and
\begin{equation} \label{dec1UUstarmart-2}
 \Vert r_{k,m}^Z\Vert_{\infty} \leq  2  \Vert \E ( |X_{k,m}  |   | {\mathcal G}_{k-1} ) \Vert_{\infty}  \leq 2   \int_G  \log (N(g)) \mu (dg ) < \infty  \, .
\end{equation}
Next, by \eqref{new3EFmr}, for any integers $k,i$ in $[2m+1, 3m]$ such that $k \geq i$, 
\begin{equation} \label{dec1UUstarmart-3}
 \Vert \E ( r_{k,m}^Z  | {\mathcal F}^Z_{i-1} ) \Vert_{\infty}  \leq  \delta_{1, \infty} (k - i+1)  \, . 
\end{equation}
In addition, for any $r \geq 1$,  $\Vert d^Z_{k,m}  \Vert_r  \leq 2 \Vert Z_{k,m} \Vert_r$ and, for any integer $k \in [2m+1, 3m]$, 
  \begin{align*}
 \E \big | {Z}_{k,m}  \big |^r 
 &  = \E \Big | f_m (\varepsilon_{k-m+1}, \ldots, \varepsilon_m, \varepsilon_{2m+1}, \ldots,  \varepsilon_{k})  - f_m (\varepsilon^*_{k-m+1}, \ldots, \varepsilon^*_{2m}, \varepsilon_{2m+1}, \ldots,  \varepsilon_{k})   \Big |^r \nonumber \\
 &  \leq   \int  \!\! \int \E \Big | f_m (u_{k-m+1}, \ldots, u_{2m}, \varepsilon_{2m+1}, \ldots,  \varepsilon_{k})  \nonumber  \\& \quad \quad  -   f_m (v_{k-m+1}, \ldots, v_{2m}, \varepsilon_{2m+1}, \ldots,  \varepsilon_{k}) \Big |^r   \prod_{i=k-m+1}^{2m} d\mu (v_i)   \prod_{i=k-m+1}^{2m} d\mu (u_i)  \nonumber \\
&  \leq \sup_{{\bar x} , {\bar y} \in X} \E  | X_{k-2m,{\bar x}} - X_{k-2m,{\bar y}}|^r= \delta^r_{r,\infty}\ (k-2m) \,  ,
\end{align*}
implying that 
\begin{equation} \label{P1lmaU1normer} 
\Vert d^Z_{k,m}  \Vert_r  \leq 2 \delta^r_{r,\infty}\ (k-2m) \, .
\end{equation}
Starting from \eqref{dec1UUstarmart}, considering the upper bound  \eqref{new1EFmr} and proceeding as in the proof of Lemma \ref{lmaR1normep} by taking into account  the upper bounds 
\eqref{dec1UUstarmart-1}-\eqref{P1lmaU1normer}, the lemma follows. \qed

  \medskip
  
 For the lemmas below, we recall the definitions \eqref{defUj}, \eqref{defRj},  \eqref{notaSm1} and \eqref{defU2*}  for $U_2$, $R_2$,  $Y_2^{(1)}$ and $U_2^*$.

 \begin{Lemma} \label{zolotarev3-step1} Let $r \in ]2,3]$. Assume that  $\mu $ has a finite moment of order $r+1$. Let  $\alpha_m = \sqrt{\frac{\E_{{\mathbb F}_m} ( (U_2+R_2)^2)}{\E_{{\mathbb F}_m} ( (U_2^{*})^2)}}$. Then for $f(x) \in \{ \cos x, \sin x \}$, we have 
\[
\Big \Vert \E_{{\mathbb F}_m}  \Big [ f \Big ( t \frac{Y_2^{(1)}}{\sqrt{2m}} \Big ) \Big ]  - \E_{{\mathbb F}_m}  \Big [ f \Big ( t \alpha_m \frac{U_2^{*}}{\sqrt{2m}} \Big ) \Big ] \Big   \Vert_{1}  \ll |t|^r m^{-1/2} \, .
\]
 \end{Lemma}
 \noindent {\bf Proof of Lemma \ref{zolotarev3-step1}.}  Note that $h = f/2^{3-r} $ is such that $| h''(x) - h''(y) | \leq |x-y|^{r-2}$. Using the arguments developed in the proof of \cite[Lemma 5.2, Item 3]{DMR09} and setting $V = U_2+R_2 - U_2^{*}$ and ${\tilde V} = V + (1 - \alpha_m) U_2^*$, we get 
 \begin{multline} \label{p1zolotarev3}
2^{r-3} (r-1) \times (2m)^{r/2} \Big | \E_{{\mathbb F}_m}  \Big [ f \Big ( t \frac{Y_2^{(1)}}{\sqrt{2m}} \Big ) \Big ]  - \E_{{\mathbb F}_m}  \Big [ f \Big ( t \alpha_m \frac{U_2^{*}}{\sqrt{2m}} \Big ) \Big ] \Big   |  \\ \leq  |t|^r \Big \{ \alpha_m^{r-1} \big (  \E_{{\mathbb F}_m}  (|{\tilde V}|^r)\big )^{1/r}
 \big (  \E (|U_2^*|^r)\big )^{(r-1)/r}  \\+ \alpha^{r-2}_m \big (  \E_{{\mathbb F}_m}  (|{\tilde V}|^r)\big )^{2/r}
 \big (  \E (|U_2^*|^r)\big )^{(r-2)/r}  + 
  \E_{{\mathbb F}_m}  (|{\tilde V}|^r) \Big \} 
  \, .
\end{multline}
Next, note that, by H\"older's inequality, 
 \begin{align*}
 \E \big ( \alpha_m^{r-1} \big (  \E_{{\mathbb F}_m}  (|{\tilde V}|^r)\big )^{1/r} \big ) & \leq   \E \big ( \alpha_m^{r-1} \big (  \E_{{\mathbb F}_m}  (|V|^r)\big )^{1/r} \big ) +    \E \big ( \alpha_m^{r-1} \times |1 - \alpha_m|  \big )  \Vert U_2^* \Vert_r  \\
&  \leq  \Vert \alpha_m \Vert_r^{r-1} \Vert V \Vert_r +  \Vert \alpha_m \Vert_r^{r-1} \Vert 1 - \alpha_m \Vert_r    \Vert U_2^* \Vert_r \, .
\end{align*}
Proceeding similarly for the two last terms in \eqref{p1zolotarev3} and taking the expectation, we derive 
 \begin{align*}
2^{r-3} (r-1) \times (2m)^{r/2}   &  \Big \Vert \E_{{\mathbb F}_m}  \Big [ f \Big ( t \frac{Y_2^{(1)}}{\sqrt{2m}} \Big ) \Big ]   - \E_{{\mathbb F}_m}  \Big [ f \Big ( t \alpha_m \frac{U^*_2}{\sqrt{2m}} \Big ) \Big ] \Big  \Vert_1  \\ &  \leq  |t|^r \Vert \alpha_m \Vert_r^{r-1} \Vert V \Vert_r   \Vert U_2^* \Vert^{r-1}_r+ |t|^r \Vert \alpha_m \Vert_r^{r-1} \Vert 1 - \alpha_m \Vert_r    \Vert U_2^* \Vert^r_r \\ &   \quad +  2 |t|^r \Vert \alpha_m \Vert_r^{r-2} \Vert V \Vert^{2}_r   \Vert U_2^* \Vert^{r-2}_r+  2  |t|^r \Vert \alpha_m \Vert^{r-2}_r \Vert 1 - \alpha_m \Vert^2_r    \Vert U_2^* \Vert^r_r \\ &  \quad +  2^{r-1} |t|^r  \Vert V \Vert^r_r  +   2^{r-1} |t|^r  \Vert 1 - \alpha_m \Vert^r_r    \Vert U_2^* \Vert^r_r
  \, .
\end{align*}
According to Lemmas \ref{lmaR1normep} and \ref{lmaU1U1starnormep}, since $\mu $ has a moment of order $r+1$, $\Vert V \Vert_r \ll 1$. Moreover, by Lemma \ref{lmamomentp}, $  \Vert U_2^* \Vert_r =  \Vert \sum_{k=1}^m X_{k+m,m}\Vert_r \leq \sqrt{m} $. On another hand,
 \begin{multline*}
 \Vert U_2^* \Vert_2 \times  \Vert 1 - \alpha_m \Vert_r   = \Big \Vert  \sqrt{\E_{{\mathbb F}_m} ( (U_2+R_2)^2) }  - \sqrt{\E_{{\mathbb F}_m} ( (U^*_2)^2) }  \Big \Vert_r \\  \leq  \Big \Vert  \sqrt{\E_{{\mathbb F}_m} ( (U_2+R_2 - U_2^*)^2) }   \Big \Vert_r  \leq \Vert V \Vert_r \ll 1 \, .
\end{multline*}
Since $ \lim_{m \rightarrow \infty} m^{-1}  \Vert U_2^* \Vert_2^2 = s^2>0$, it follows that for $m$ large enough 
\beq \label{diff1alpha} \Vert 1 - \alpha_m \Vert_r  \ll m^{-1/2}  \, .
\eeq The lemma follows from all the above considerations.  \qed

 \begin{Lemma} \label{zolotarev3-step2} Let $r \in ]2,3]$. Assume that  $\mu $ has a finite moment of order $q=r+1$. Recall the notation  $\alpha_m = \sqrt{\frac{\E_{{\mathbb F}_m} ( (U_2+R_2)^2)}{\E_{{\mathbb F}_m} ( (U^*_2)^2)}}$. Then for $f(x) \in \{ \cos x, \sin x \}$, we have 
\[
\Big \Vert \E_{{\mathbb F}_m}  \Big [ f \Big ( t \alpha_m \frac{U^*_2}{\sqrt{2m}} \Big ) \Big ]  -  \E_{{\mathbb F}_m}  \Big [  f \Big ( t  \alpha_m  \frac{s_m N}{\sqrt{2}}  \Big )  \Big ]  \Big   \Vert_{1}  \ll |t|^r m^{-1/2} + |t| m^{-(r-1/2)} \, ,
\]
where  $s_m^2 = \E (S_m^2) /m$ and $N $ is a standard Gaussian random variable independent of ${\mathbb F}_m$. 
 \end{Lemma}
 \noindent {\bf Proof of Lemma \ref{zolotarev3-step2}.} Let $W_0^*$ be distributed as $W_0$ and independent of $W_0$. Let $(\varepsilon^*_k)_{k \geq 1}$ be an independent copy of $(\varepsilon_k)_{k \geq 1}$, independent of $(W_0^* ,W_0)$. Define $S^*_m = \sum_{k=m+1}^{2m} X_k^*$ where $X_k^*= \sigma ( \varepsilon_k^*, W^*_{k-1}) - \lambda_{\mu}$ with $W^*_{k}=\varepsilon_k^* W^*_{k-1}$, for $k \geq 1$. Note that  $S^*_m$ is independent of 
 ${\mathbb F}_m$ and has the same law as $S_m$. In addition,  by stationarity,  \cite[Lemma 24]{CDM} (applied with  $M_k = +\infty$) and \eqref{estimatedeltabis}, 
  \begin{multline} \label{zolotarev3-step2-p1}
 \Big \Vert \E_{{\mathbb F}_m}  \Big [ f \Big ( t \alpha_m \frac{S_m^{*}}{\sqrt{2m}} \Big ) \Big ]   - \E_{{\mathbb F}_m}  \Big [ f \Big ( t \alpha_m \frac{U^*_2}{\sqrt{2m}} \Big ) \Big ]  \Big ]  \Big   \Vert_{1}  \\
   \ll 
\frac{|t|}{\sqrt{2m}}  \E |\alpha_m | \times   \sum_{k=m+1}^{2m}  \Vert  X_{k,m} -X_k \Vert_1   \ll  \frac{|t|}{\sqrt{m}} \times m \delta_{1, \infty} (m) \ll  |t| m^{-(r-1/2) } \, .
\end{multline}
On another hand, let  $h = f/2^{3-r} $ and note that  $| h''(x) - h''(y) | \leq |x-y|^{r-2}$.  Hence, by the definition of the Zolotarev distance of order $r$, 
\[
 \Big \Vert \E_{{\mathbb F}_m}  \Big [ f \Big ( t \alpha_m \frac{S_m^{*}}{\sqrt{2m}} \Big ) \Big ]  -  \E_{{\mathbb F}_m}  \Big [  f \Big ( t  \alpha_m  \frac{s_m N}{\sqrt{2}}  \Big )  \Big ] \Big   \Vert_{1}  \leq  2^{3-r} |t|^r \times \Vert \alpha_m \Vert_r^r \zeta_r \big (  P_{S_m/{\sqrt{ 2  m}}} , G_{s_m^2/2} \big) \, .
\]
Next we apply \cite[Theorem 3.2, Item 3.]{DMR09} and derive that since $\mu$ has a moment of order $q>3$, 
\[
 \zeta_r \big (  P_{S_m/{\sqrt{2 m}}} , G_{s_m^2/2} \big) \ll m^{-1/2} \, .
\]
As we mentioned before, the fact that the conditions (3.1), (3.4) and (3.5) required in  \cite[Theorem 3.2]{DMR09} hold when $\mu$ has a moment of order  $q>3$ has been proved in the proof of  \cite[Theorem 2]{CDJ}. 
Hence, since  $ \Vert \alpha_m \Vert_r \ll 1$ (see \eqref{diff1alpha}), 
\begin{equation} \label{zolotarev3-step2-p2}
 \Big \Vert \E_{{\mathbb F}_m}  \Big [ f \Big ( t \alpha_m \frac{S_m^{*}}{\sqrt{2m}} \Big ) \Big ]  -  \E_{{\mathbb F}_m}  \Big [  f \Big ( t  \alpha_m  \frac{s_m N}{\sqrt{2}}  \Big )  \Big ] \Big   \Vert_{1}  \ll  \frac{| t|^r }{\sqrt{m}}    \, .
\end{equation}
Considering the upper bounds  \eqref{zolotarev3-step2-p1} and   \eqref{zolotarev3-step2-p2}, the lemma follows. \qed

 \begin{Lemma} \label{zolotarev3-step3} Let $r \in ]2,3]$. Assume that  $\mu $ has a finite moment of order $q=r+1$. Recall the notations  $\alpha_m = \sqrt{\frac{\E_{{\mathbb F}_m} ( (U_2+R_2)^2)}{\E_{{\mathbb F}_m} ( (U^*_2)^2)}}$ and $s_m^2 = \E (S_m^2) /m$. . Then, for $f(x) \in \{ \cos x, \sin x \}$,  
\[
\Big \Vert  \E_{{\mathbb F}_m}  \Big [  f \Big ( t  \alpha_m  \frac{s_m N}{\sqrt{2}}  \Big )  \Big ]  -  \E_{{\mathbb F}_m}  \Big [  f \Big ( t   \frac{s N}{\sqrt{2}}  \Big )  \Big ]  \Big   \Vert_{1}  \ll  \frac{ |t| }{ m^{1/2 + \eta } }  \, .
\]
where $\eta = \min ( \frac{3}{10} ,  \frac{r-2}{2} , \frac{r-2}{2r-3} )$ and $N $ is a standard Gaussian random variable independent of ${\mathbb F}_m$. 
 \end{Lemma}
 
  \noindent {\bf Proof of Lemma \ref{zolotarev3-step3}.} We  have
  \begin{multline} \label{zolotarev3-step3P1}
  \Big \Vert  \E_{{\mathbb F}_m}  \Big [  f \Big ( t  \alpha_m  \frac{s_m N}{\sqrt{2}}  \Big )  \Big ]  -  \E_{{\mathbb F}_m}  \Big [  f \Big ( t   \frac{s N}{\sqrt{2}}  \Big )  \Big ]  \Big   \Vert_{1}  \\
   \leq |t |    \E |N|  \big (   \Vert \alpha_m \Vert_1 | s - s_m | + s  \times\Vert 1 - \alpha_m \Vert_1 \big )   \, .
\end{multline}
 But, since  $ \lim_{m \rightarrow \infty} m^{-1}  \Vert U_2^* \Vert_2^2 = s^2>0$, 
 \[
\Vert 1 - \alpha_m \Vert_1  \leq   \Vert 1 - \alpha^2_m \Vert_1 \sim \frac{1}{s^2 m} \big \Vert  \E_{{\mathbb F}_m} ( (U_2+R_2)^2) - \E_{{\mathbb F}_m} ( (U_2^*)^2)\big \Vert_1 \, .
 \]
On another hand 
  \[
\big \Vert  \E_{{\mathbb F}_m} ( (U_2+R_2)^2) - \E_{{\mathbb F}_m} ( (U_2^*)^2)\big \Vert_1 \leq  \big \Vert  \E_{{\mathbb F}_m}   (U_2^2) - \E ( (U_2^*)^2) \big \Vert_1 + \Vert R_2 \Vert_2^2 + 2 \Vert   \E_{{\mathbb F}_m} ( U_2R_2) \Vert_1 \, .
 \]
But, by stationarity, 
  \begin{multline*}
  \big \Vert  \E_{{\mathbb F}_m}   (U_2^2) -  \E ( (U_2^*)^2)  \big \Vert_1 \leq  \Big \Vert \E_m \Big ( \sum_{k=m+1}^{2m} { \bar X}_{k,m}  \Big )^2  -  \E  \Big ( \sum_{k=m+1}^{2m} { \bar X}_{k,m}  \Big )^2  \Big \Vert_{1}  \\+ \Big ( \sum_{k=2m+1}^{3m} \Vert  \E_{{\mathbb F}_m} (X_{k,m} ) \Vert_2 \Big )^2 \, ,
  \end{multline*}
  where ${ \bar X}_{k,m} = X_{k,m} - \E_m (X_{k,m} )$ and $\E_m (\cdot) = \E(\cdot |   \sigma (W_0, \varepsilon_1, \ldots, \varepsilon_m))$. 
 Hence, by \eqref{new1EFmr} and Lemma \ref{lma4.7}, since $q=r+1 $ and $r >2$,  
 \[
 \big \Vert  \E_{{\mathbb F}_m}  (U_2^2) - \E( U_2^2) \big \Vert_1  \ll m^{1/5}\, .
 \]
By stationarity and Lemma \ref{lmaR1normep}, we also have $\Vert R_2 \Vert_2= \Vert R_1 \Vert_2 \ll 1$. Therefore
  \begin{equation*} \label{zolotarev3-step3P1bis}
\big \Vert  \E_{{\mathbb F}_m} ( (U_2+R_2)^2) - \E_{{\mathbb F}_m} ( (U_2^*)^2)\big \Vert_1 \ll  m^{1/5} +  \Vert   \E_{{\mathbb F}_m} ( U_2R_2) \Vert_1 \, .
 \end{equation*}
Next, note that 
\[
 \Vert   \E_{{\mathbb F}_m} ( U_2R_2) \Vert_1 = \Big    \Vert   \E_{{\mathbb F}_m}  \Big (  R_2 \sum_{k=2m+1}^{3m} X_{k,m}\Big )  \Big \Vert_1  \, .
\]
Let $h(m) $ be a positive  integer less than $m$. Using stationarity, Lemma \ref{lmaR1normep} and similar arguments as those developed in the proof of Lemma \ref{lmamomentp},  we first notice that  
\[
\Big    \Vert   \E_{{\mathbb F}_m}  \Big (  R_2 \sum_{k=3m-h(m) +1}^{3m} X_{k,m}\Big )  \Big \Vert_1  \leq \Vert R_2 \Vert_2 \Big \Vert  \sum_{k=3m-h(m) +1}^{3m} X_{k,m}  \Big \Vert_2 \ll \sqrt{h(m) } \, .
\]
We handle now the term $  \Vert   \E_{{\mathbb F}_m}  \big (  R_2 \sum_{k= 2m+1}^{3m-h(m) }X_{k,m}\big )  \Vert_1$. 
 For $2m+1 \leq k \leq 3m$, define $X_{k,m}^* $ as in \eqref{defU2*}. Using \eqref{new2EFmr}  and \eqref{estimatedeltabis}, note that 
\[
 \sum_{k= 2m +1}^{3m - h(m) }   \Vert  X_{k,m} -X_{k,m}^*  \Vert_{2} \leq  \sum_{k=2m+1}^{3m} \delta_{2, \infty} (k-2m)  \ll  \sum_{k=1}^{m} k^{-(q/2-1)}  \, . \]
 Hence
 \[
 \sum_{k= 2m +1}^{3m - h(m) }   \Vert  X_{k,m} -X_{k,m}^*  \Vert_{2}  \ll  m^{(3-r) /2}    {\bf 1}_{r < 3} +  {\bf 1}_{r = 3}  \log (m)   \, .
 \]
This estimate combined with $\Vert R_2 \Vert_2 \ll 1$ entails
\[
\Big    \Vert   \E_{{\mathbb F}_m}  \Big (  R_2  \sum_{k= 2m +1}^{3m - h(m) } X_{k,m}\Big )  \Big \Vert_1  \ll m^{(3-r) /2}    {\bf 1}_{r < 3} +  {\bf 1}_{r = 3}  \log (m) +  \Big    \Vert   \E_{{\mathbb F}_m}  \Big (  R_2  \sum_{k= 2m +1}^{3m - h(m) } X^*_{k,m}\Big )  \Big \Vert_1 .
\]
Since $(X_{k,m}^* )_{2m+1 \leq k \leq 3m}$ is independent of ${\mathbb F}_m$, we have  $\E ( X_{k,m}^* | {\mathbb F_m} )=0$ for any $2m+1 \leq k \leq 3m$. Hence 
\[ \Big    \Vert   \E_{{\mathbb F}_m}  \Big (  R_2  \sum_{k= 2m +1}^{3m - h(m) } X^*_{k,m}\Big )  \Big \Vert_1 =   \Big    \Vert   \E_{{\mathbb F}_m}  \Big (  \sum_{k=2m+1}^{3m-h(m)} X^*_{k,m}  \sum_{\ell=3m +  1}^{4m} X_{\ell,m}  \Big )  \Big \Vert_1 \, .
\]
Next, note that if $\ell-m+1 \geq k+1$, conditionally to  ${\mathbb F}_m$, $X^*_{k,m} $ is independent of $ X_{\ell,m} $, which implies that $  \E_{{\mathbb F}_m}(  X^*_{k,m} X_{\ell,m} ) =0$.  Hence
\[
 \Big    \Vert   \E_{{\mathbb F}_m}  \Big (  \sum_{k=2m+1}^{3m-h(m) } X^*_{k,m}  \sum_{\ell=3m+1}^{4m} X_{\ell,m}  \Big )  \Big \Vert_1 = \Big    \Vert   \E_{{\mathbb F}_m}  \Big (  \sum_{k=2m+1}^{3m-h(m)} X^*_{k,m}  \sum_{\ell=3m+1}^{4m-h(m) -1} X_{\ell,m}  \Big )  \Big \Vert_1 \, .
\]
Now, for any $3m+1 \leq \ell \leq 4m - h(m) -1$, let 
\[
X_{\ell,m}^{ (h(m) ,*)} = f_m ( \varepsilon^*_{\ell-m+1}, \ldots, \varepsilon^*_{3m-h(m)}, \varepsilon_{3m-h(m) +1}, \ldots \varepsilon_\ell  ) \, ,
\]
and note that $ \E_{{\mathbb F}_m} ( X^*_{k,m} X_{\ell,m}^{ (h(m),*)}  ) =0$ for any $k \leq 3m-h(m)$ and any $\ell \geq 3m+1$.  So, overall,  setting $q' = q/(q-1)$, 
 \begin{multline*}
 \Big    \Vert   \E_{{\mathbb F}_m}  \Big (  R_2 \sum_{k=2m+1}^{3m-h(m)} X^*_{k,m}\Big )  \Big \Vert_1  = \Big    \Vert   \E_{{\mathbb F}_m}  \Big (  \sum_{k=2m+1}^{3m-h(m) } X^*_{k,m}  \sum_{\ell=3m+1}^{4m - h(m) -1} ( X_{\ell,m}  -  X_{\ell,m}^{ (h(m) ,*)}  ) \Big )  \Big \Vert_1 \\
 \leq  \Big \Vert  \sum_{k=2m +1}^{3m-h(m) } X^*_{k,m}  \Big \Vert_q    \sum_{\ell=3m+1}^{4m - h(m) -1} \Vert X_{\ell,m}  -  X_{\ell,m}^{ (h(m),*)}   \Vert_{q'} \, .
 \end{multline*}
Proceeding as in the proof of \cite[Lemma 24]{CDM}, we infer that the following inequality holds: $\Vert X_{\ell,m}  -  X_{\ell,m}^{ (h(m),*)}   \Vert_{q'}  \leq  \delta_{q', \infty} (\ell  - 3m + h(m))$.  
 Hence, taking into account \eqref{estimatedeltabis} and Lemma \ref{lmamomentp}, we get 
 \[
 \Big    \Vert   \E_{{\mathbb F}_m}  \Big (  R_2 \sum_{k=2m+1}^{3m-h(m) } X^*_{k,m}\Big )  \Big \Vert_1  \ll  \sqrt{m} \sum_{\ell \geq h(m)}\frac{1}{\ell^{q-2}} \ll  \sqrt{m} (h(m))^{2-r} \, .
 \]
Taking into account  all the above considerations and selecting $h(m) = m^{1/(2r-3)}$, we derive 
  \begin{equation} \label{zolotarev3-step3P2}
m   \Vert 1 - \alpha_m \Vert_1    \ll     m^{(3-r)/2}    {\bf 1}_{r < 3}   + m^{1/ (4r-6)}  + m^{1/5}  \,  .
  \end{equation}
  On another hand, since $s^2 >0$, $ |   s-s_m | \leq s ^{-1}  |   s^2-s^2_m |$. Hence by using Remark \ref{remonvar}, the definition of $s^2_m$ and stationarity, we derive that 
  \[
|   s-s_m | \leq  \frac{2  }{s m }  \sum_{k\geq 1} k |{\rm Cov} (X_0,X_k)|  \, .
  \]
By the definition of $\delta_{1, \infty}(k)$, $ |{\rm Cov} (X_0,X_k)|  \leq \Vert X_0 \Vert_1 \delta_{1, \infty} (k) $. So, by using \eqref{estimatedelta} and the fact that $q \geq 2$, we get
   \begin{equation} \label{zolotarev3-step3P3}
 |   s-s_m |   \ll  m^{-1} \,  .
  \end{equation}
  Starting from \eqref{zolotarev3-step3P1} and taking into account  \eqref{zolotarev3-step3P2} and  \eqref{zolotarev3-step3P3}, the lemma follows. \qed

\medskip

Combining Lemmas \ref{zolotarev3-step1}, \ref{zolotarev3-step2} and \ref{zolotarev3-step3}, we derive 
 \begin{Lemma} \label{lma4.9q=4} Let $r \in ]2,3]$. Assume that  $\mu $ has a finite moment of order $q=r+1$.  Then, for $f(x) \in \{ \cos x, \sin x \}$,
 \[
\Big \Vert \E_{{\mathbb F}_m}  \Big [ f \Big ( t \frac{Y_2^{(1)}}{\sqrt{2m}} \Big ) \Big ]  -  \E \big [ f  ( t  s N /{\sqrt 2} ) \big ] \Big   \Vert_{1}  \ll |t|^r m^{-1/2} + |t| m^{-(1/2 + \eta )} \,,
\]
where $\eta = \min ( \frac{3}{10} ,  \frac{r-2}{2} , \frac{r-2}{2r-3} )$.
\end{Lemma}
Let $R_1$ be defined by \eqref{defRj}. Proceeding similarly as to derive the previous lemma,  we get 
 \begin{Lemma} \label{lma4.9q=4bis}Let $r \in ]2,3]$. Assume that  $\mu $ has a finite moment of order $q=r+1$.  Then for $f(x) \in \{ \cos x, \sin x \}$,  
\[
\Big \Vert \E_{{\mathbb F}_m}  \Big [ f \Big ( t \frac{ \sum_{k=1}^m  X_{k}   +R_1}{\sqrt{2m}} \Big ) \Big ]  -  \E \big [ f  ( t  s N /{\sqrt 2} ) \big ] \Big   \Vert_{1}  \ll |t|^r m^{-1/2} + |t| m^{-(1/2 + \eta )} \, ,
\]
where $\eta = \min ( \frac{3}{10} ,  \frac{r-2}{2} , \frac{r-2}{2r-3} )$.
\end{Lemma}

\section*{Acknowledgements} 

This research was partially supported by the NSF grant DMS-2054598. The authors would like to thank two anonymous referees for their valuable suggestions, which improved the presentation of the paper.


\begin{thebibliography}{4}




\bibitem{BQ} Benoist, Y. and Quint, J.-F. (2016). Central limit theorem for linear groups, {\em Ann. Probab.} 
{\bf 44} no. 2, 1308--1340.

\bibitem{Bo82} Bolthausen, E. (1982). 
Exact convergence rates in some martingale central limit theorems.
\textit{Ann. Probab.} \textbf{ 10}, no. 3, 672--688. 

\bibitem{BL} Bougerol, P. and  Lacroix, J.  Products of random matrices with applications to Schr\"odinger operators. Progress in Probability and Statistics, 8. Birkh\"auser Boston, Inc., Boston, MA,1985.

\bibitem{Bu73}  Burkholder, D. L. (1973). Distribution function inequalities for martingales. \textit{Ann. Probab.} \textbf{ 1}, 19--42.

\bibitem{CDJ} Cuny, C., Dedecker, J. and Jan, C. (2017). Limit theorems for the left random walk on $GL_d({\mathbb R})$.  \textit{Ann. Inst. H. Poincar\'e Probab. Statist.} \textbf{53}, no. 4, 1839--1865.

\bibitem{CDM} Cuny, C., Dedecker, J. and Merlev\`ede, F. (2018). On the Koml\'os, Major and Tusn\'ady strong approximation for some classes of random iterates.  \textit{Stochastic Process. Appl.}  \textbf{128}, no. 4, 1347--1385.

\bibitem{CDMP} Cuny, C., Dedecker, J.,  Merlev\`ede, F. and Peligrad, M. (2022). Berry-Esseen type bounds for the matrix coefficients and the spectral radius of the left random walk on $GL_d({\mathbb R})$. 
 \textit{C. R. Math. Acad. Sci. Paris}  \textbf{360}, 475--482.


%\bibitem{De10}  Dedecker, J. (2010). An empirical central limit theorem for intermittent maps. \textit{Probab. Theory Related Fields} \textbf{ 148},   no. 1-2, 177--195.

\bibitem{DD03}  Dedecker, J. and Doukhan, P. (2003). A new covariance inequality and applications. \textit{Stochastic Process. Appl.}  \textbf{106}, no. 1, 63--80.

\bibitem{DMR09} Dedecker, J.,  Merlev\`ede, F. and Rio, E. (2009). Rates of convergence for minimal distances
in the central limit theorem under projective criteria. \textit{Electron. J. Probab.} \textbf{14}, no. 35, 978--1011.

\bibitem{DKW} Dinh, T.-C., Kaufmann, L. and Wu, H. Random walks on $SL_2({\mathbb C})$: spectral gap and local
limit theorems. https://arxiv.org/pdf/2106.04019.pdf

\bibitem{FP20} Fernando, K. and P\`ene, F.  (2022). Expansions in the local and the central limit theorems for dynamical systems. \textit{Commun. Math. Phys.} \textbf{389}, 273--347.



\bibitem{GR} Guivarc'h, Y.  and  Raugi, A. (1985). Fronti\`ere  de Furstenberg,  propri\'et\'es  de contraction  et th\'eor\`emes de convergence, {\em Z. Wahrsch. Verw. Gebiete} {\bf 69} no. 2, 187-242.

%\bibitem{HH80} Hall, P. and  Heyde, C. C.
%Martingale limit theory and its application.
%Probability and Mathematical Statistics. Academic Press, Inc., New York-London, 1980. 

%\bibitem{Ho} H\"ormann, S.  (2009). Berry–Esseen bounds for econometric time series. \textit{ALEA Lat. Am. J.
%Probab. Math. Stat.} \textbf{ 6} 377--397.

\bibitem{Fe71} Feller, W. 
{\it An introduction to probability theory and its applications}. Vol. II.
Second edition John Wiley $\&$ Sons, Inc., New York-London-Sydney 1971 xxiv+669 pp. 

\bibitem{FK} Furstenberg, H. and Kesten, H.  (1960). Products of Random Matrices. 
{\it Ann. Math. Statist.} {\bf 31}, no. 2, 457--469.

\bibitem{Hoeffding} Hoeffding, W. (1963). Probability inequalities for sums of bounded random variables. \textit{J. Amer. Statist. Assoc.} \textbf{58}, 13--30.

\bibitem{Jan} Jan,  C. (2001). Vitesse de convergence dans le TCL pour des processus associ\'es \`a des syst\`emes dynamiques ou  des  produits  de  matrices  al\'eatoires,  Th\`ese de l'Universit\'e de Rennes 1 (2001),  thesis number 01REN10073

\bibitem{Ji}  Jirak, M. (2016). Berry-Esseen theorems under weak dependence. {\it Ann. Probab.} {\bf  44}, no. 3, 2024--2063.

\bibitem{Ji20}  Jirak, M. (2020).  A Berry-Esseen bound with (almost) sharp dependence conditions. \textit{arXiv}:1606.01617

\bibitem{LP} Le Page, E. (1982). Th\'eor\`emes  limites  pour  les  produits  de  matrices  al\'eatoires, Probability measures on groups (Oberwolfach, 1981), pp. 258--303, Lecture Notes in Math., 928, Springer, Berlin-New York.

\bibitem{MPU19}  Merlev\`ede, F., Peligrad, M. and Utev, S. \textit{Functional Gaussian approximation for dependent structures}. Oxford Studies in Probability, 6. Oxford University Press, Oxford, 2019. xv+478 pp

\bibitem{XGLeuropean} Xiao,  H. Grama, I. and Liu, Q. (2021). Berry-Esseen bound and precise moderate deviations
for products of random matrices.  \textit{Journal of the European Mathematical Society}, European Mathematical
Society, In press, 10.4171/JEMS/1142. hal-03431385.

\bibitem{XGL} Xiao,  H. Grama, I. and Liu, Q. (2021). Berry Esseen bounds and moderate deviations for  random walks on $GL_d({\mathbb R})$.  \textit{Stochastic Process. Appl.}  \textbf{142}, 293--318.

\end{thebibliography}
\end{document}